\documentclass[11pt]{article}
\usepackage{amssymb,amsmath}
\usepackage{graphics}
\usepackage{graphicx}
\def\K{\mathbb{K}}

\def\1{\mbox{\bf 1}} 

\newcommand\Lip{{\mathrm {Lip } \, }}

\newtheorem{theo}{Theorem}
\newtheorem{Lemma}{Lemma}
\newtheorem{Ex}{Example}
\newtheorem{Rem}{Remark}
\newtheorem{Def}{Definition}
\newtheorem{proofl}{Proof of the lemma}

\def\1{\mbox{I\hspace{-.6em}1}} 

\newenvironment{hyp}[1]{\refstepcounter{hyp}\begin{itemize}\item[{\bf
      (A\arabic{hyp})}] \label{hyp:#1}}{\end{itemize}}
\newcounter{hyp}

\author{Paul Doukhan, Silika Prohl, Christian Y. Robert}
\date{\today}
\title{Subsampling weakly dependent times series\\ and application to extremes}

\begin{document}
\date{\today}
\maketitle

\begin{abstract}
This paper provides extensions of the work on subsampling by Bertail \textit{%
et al.} (2004) for strongly mixing case to weakly dependent case by
application of the results of Doukhan and Louhichi (1999). We
investigate properties of smooth and rough subsampling estimators
for distributions of converging and extreme statistics when the
underlying time series is $\eta $ or $\lambda $-weakly dependent.
\end{abstract}

\section{Introduction}

Politis and Romano (1994) \cite{PR} established the subsampling estimator
for statistics when the underlying sequence is strongly mixing. Bertail 
\textit{et al.} (2004) \cite{BHPW} applied this work to subsampling
estimators for distributions of diverging statistics. In particular, they
constructed an approximation of the distribution of the sample maximum
without any information on the tail of the stationary distribution. However
the assumption on the strong mixing properties of the time series is
sometimes too strong as for the class of first-order autoregressive
sequences introduced and studied by Chernick (1981) \cite{Ch}: for $t\in 
\mathbb{Z}$, let $X_{t}$ be given by 
\begin{equation}
X_{t}=\frac{1}{r}(X_{t-1}+\varepsilon _{t}),  \label{FOS}
\end{equation}
where $r\geq 2$ is an integer, $(\varepsilon _{t})_{t\in \mathbb{Z}}$ are
iid and uniformly distributed on the set $\{0,1,\ldots ,r-1\}$ and $X_{0}$
is uniformly distributed on $[0,1]$. Andrews (1984) \cite{And} and Ango-Nze
and Doukhan (2004) \cite{AD} (see page 1009 and Note 5 on page 1028) give
arguments to derive that such models are not mixing. The results of Bertail 
\textit{et al.} (2004) \cite{BHPW} can not be used although the normalized
sample maximum has a non degenerate limiting distribution: let $M_{n}=\max
(X_{1},\ldots ,X_{n})$, then 
\begin{equation*}
\lim_{n\rightarrow \infty }\mathbb{P}\big(n(1-M_{n})\leq x\big)=1-\exp \big(%
-r^{-1}(r-1)x\big),\qquad \text{for all }x\geq 0\text{,}
\end{equation*}
(see Theorem 4.1 in Chernick (1981) \cite{Ch}).

This paper is aimed at weakening the dependence conditions assumed in
Bertail \textit{et al.} (2004) \cite{BHPW} and at studying new smooth
subsampling estimators adapted to our weak dependence conditions.

Doukhan and Louhichi (1999) \cite{dl} introduced a wide dependence framework
that turns out in particular to apply to the previous processes and that
widely improves the amount of potentially usable models. This dependence
structure is addressed in Section \ref{wd}. In Section \ref{subsam} we
introduce smooth and rough subsampling estimators for the distribution of
converging statistics and studied their asymptotic properties. We consider
two subsampling schemes based on overlapping and non-overlapping samples. In
the next section we consider subsampling estimators for the distribution of
extremes and, to fix ideas, we focus on the case of the normalized sample
maximum. We first discuss sufficient conditions adapted to our weak
dependence framework such that the normalized maximum converges in
distribution. Then we discuss how to estimate the normalizing sequences and
we derive the asymptotic properties of the subsampling estimators. A
simulation study provides explicit comparisons of the various considered
subsamplers in Section \ref{SEx}. Proofs are reported in a last section.

\section{Weak dependence}

\label{wd}

Doukhan and Louhichi (1999) \cite{dl} proposed a new idea of weak dependence
that makes explicit the asymptotic independence between past and future. Let
us consider a strictly stationary time series $X=(X_{t})_{t\in \mathbb{Z}}$
which (for simplicity) will be assumed to be real-valued. Let us denote by $%
F $ its stationary distribution function. If $X$ is a sequence of iid random
variables, then for all $t_{1}\neq t_{2}$, independence between $X_{t_{1}}$
and $X_{t_{2}}$ writes $\mbox{Cov}\,(f(X_{t_{1}}),g(X_{t_{2}}))=0$ for all $%
f,g$ with $\| f\| _{\infty },\| g\| _{\infty }\leq 1$, where $\| f\|
_{\infty }$ denotes the surpremum norm of $f$. For a sequence of dependent
random variables, we would like that $\mbox{Cov}\,(f($`$past$'$),g($`$future$%
'$))$ is small when the distance between the past and the future is
sufficiently large. \newline
More precisely, for $h:\mathbb{R}^{u}\rightarrow \mathbb{R}$ ($u\in \mathbb{N%
}^{\ast }$) we define

\centerline{$\displaystyle \Lip
h=\sup_{(y_1,\ldots,y_u)\ne(x_1,\ldots,x_u) \in \mathbb{R}^u }
\frac{\left| h(y_1,\ldots,y_u)- h(x_1,\ldots,x_u) \right|}
{\|y_1-x_1\|+\cdots+\|y_u-x_u\|}.$}

\begin{Def}
\textbf{\cite{dl}}\label{defdepf} The process $X$ is $(\varepsilon ,\Psi )$%
-weakly dependent process if, for some classes of functions $\mathcal{F}_{u},%
\mathcal{G}_{v}$, $E^{u},E^{v}\rightarrow \mathbb{R},u,v\geq 1$: 
\begin{equation*}
\varepsilon (r)=\sup \frac{\left| {{\mbox{Cov}\,}}\Big(f(X_{s_{1}},\ldots
,X_{s_{u}}),g(X_{t_{1}},\ldots ,X_{t_{v}})\Big )\right| }{\Psi (f,g)}%
\rightarrow _{r\rightarrow \infty }0
\end{equation*}
where the $\sup $ bound is relative to $u,v\geq 1$, $s_{1}\leq \cdots \leq
s_{u}\leq t_{1}\leq \cdots \leq t_{v}$ with $r=t_{1}-s_{u}$, and $f,g$
satisfy $\mbox{Lip}\,f,\mbox{Lip}\,g<\infty $ and $\Vert f\Vert _{\infty
}\leq 1$, $\Vert g\Vert _{\infty }\leq 1$.

The following distinct functions $\Psi $ yield $\eta $, and $\lambda $ weak
dependence coefficients:

$
\begin{array}{lllcr}
\text{if }\Psi (f,g) & = & u\mbox{Lip}\,f+v\mbox{Lip}\,g, & \mbox{then} & 
\epsilon (r)=\eta (r), \\ 
& = & u\mbox{Lip}\,f+v\mbox{Lip}\,g+uv\mbox{Lip}\,f\cdot \mbox{Lip}\,g, & %
\mbox{then} & \epsilon (r)=\lambda (r),
\end{array}
$
\end{Def}

Note that $\lambda $-weak dependence includes $\eta $-weak dependence. A
main feature of Definition \ref{defdepf} is to incorporate a much wider
range of classes of models than those that might be described through a
mixing condition (\textit{i.e.} $\alpha $-mixing, $\beta $-mixing, $\rho $%
-mixing, $\phi $-mixing, \ldots, see Doukhan (1994) \cite{D})) or
association condition (see Chapters 1-3 in Dedecker \textit{et al.} (2007) 
\cite{df}). Limit theorems and very sharp results have been proved for this
class of processes (see Chapters 6-12 in Dedecker \textit{et al.} (2007) 
\cite{df} for more information).

We now provide a non-exhaustive list of weakly dependent sequences with
their weak dependence properties. This will prove how wide is the range of
potential applications.

\begin{Ex}
\begin{itemize}
\item  The Bernoulli shift with independent inputs $(\xi _{t})_{t\in \mathbb{%
Z}}$ is defined as $X_{t}=H((\xi _{t-j})_{j\in \mathbb{Z}}),$ $H:\mathbb{R}^{%
\mathbb{Z}}\rightarrow \mathbb{R},$ $(\xi _{i})_{i\in \mathbb{Z}}$ iid. The
process $(X_{t})_{t\in \mathbb{Z}}$ is $\eta $-weakly dependent with $\eta
(r)=2\delta _{\lbrack r/2]}^{m\wedge 1}$ if 
\begin{equation*}
\mathbb{E}\big|H\big(\xi _{j},j\in \mathbb{Z}\big)-H\big(\xi _{j}%
\mbox{I\hspace{-.6em}1}_{|j|<r},j\in \mathbb{Z}\big)\big|\leq \delta
_{r}\downarrow 0\quad (r\uparrow \infty ).
\end{equation*}
Two particular (causal) examples are given by:

- The first-order autoregressive sequences with discrete innovations given
by (\ref{FOS}). This process is not strongly mixing but it is $\eta $-weakly
dependent process such that $\eta (k)=O\left( r^{-k}\right) $.

- The LARCH model with Rademacher iid inputs: 
\begin{equation}
X_{t}=\xi _{t}(1+aX_{t-1}),\qquad \mathbb{P}(\xi _{0}=\pm 1)=\frac{1}{2}.
\label{Larch}
\end{equation}
If $a<1$, there exists a unique stationary solution (see Dedecker \textit{et
al.} (2007) \cite{df}). Doukhan, Mayo and Truquet (2008) \cite{DMT} proved
that if $a\in \left( (3-\sqrt{5})/2,1/2\right] $ the stationary solution $%
X_{t}=\xi _{t}+\sum_{j\geq 1}a^{j}\xi _{t}\cdots \xi _{t-j}$ is not strongly
mixing, but $X$ is a $\eta $-weakly dependent process such that $\eta
(k)=O\left( a^{k}\right) $.

\item  If $X$ is either a Gaussian or an associated process, then $X$ is $%
\lambda $-weakly dependent and 
\begin{equation*}
\lambda (r)=O\left( \sup_{i\geq r}|\mbox{Cov}\,(X_{0},X_{i})|\right) 
\end{equation*}
(see Doukhan and Louhichi (1999) \cite{dl}).

\item  If $X$ is a $GARCH(p,q)$ process or, more generally, a $ARCH(\infty )$
process such that $X_{k}=\rho _{k}\xi _{k}$ with $\rho
_{k}^{2}=b_{0}+\sum_{j=1}^{\infty }b_{j}X_{k-j}^{2}$ for $k\in \mathbb{Z}$
and if,

\qquad - it exists $C>0$ and $\mu \in ]0,1[$ such that $\forall j\in \mathbb{%
N}$, $0\leq b_{j}\leq C\mu ^{j}$, then $X$ is $\lambda $-weakly dependent
process with $\lambda (r)=O\left( e^{-c\sqrt{r}}\right) $ and $c>0$ (this is
the case of $GARCH(p,q)$ processes).

\qquad - it exists $C>0$ and $\nu >1$ such that $\forall j\in \mathbb{N}$, $%
0\leq b_{j}\leq Cj^{-\nu }$, then $X$ is $\lambda $-weakly dependent process
with $\lambda (r)=O\left( e^{-r}\right) $.
\end{itemize}
\end{Ex}

\section{Subsampling the distribution of converging statistics}

\label{subsam} Politis and Romano (1994) \cite{PR} introduced the
methodology of \textquotedblleft subsampling\textquotedblright\ to give
consistent approximations of confidence intervals for some parameters of the
distribution of the observations. They established the validity of their
methodology for general strongly mixing sequences under the assumption that
the considered statistics converge with a known rate.

We consider here a sequence of statistics $S_{n}=s_{n}(X_{1},\ldots ,X_{n})$
for $n=1,2,\ldots $ Let $\mathbb{K}_{n}$ be the cumulative distribution
function of $S_{n}$, $\mathbb{K}_{n}(x)=P(S_{n}\leq x)$. We assume that $%
S_{n}$ is a sequence of converging statistics in the sense that $\mathbb{K}%
_{n}$ has a limit which is denoted by $\mathbb{K}$. We assume that the
statistics satisfy one of the two following assumptions:

\strut

\noindent - \textbf{Convergent statistics:} 
\begin{eqnarray}
r_{n} &=&\sup_{x\in \mathbb{R}}\Big|\mathbb{K}_{n}(x)-\mathbb{K}(x)\Big|%
\rightarrow _{n\rightarrow \infty }0,\qquad \| \mathbb{K}^{\prime
}\|_{\infty }<\infty .  \label{con} \\
&&\qquad\quad%
\mbox{where $\K'$ denotes the density of this limit
distribution. }  \notag
\end{eqnarray}

\noindent - \textbf{Concentration condition:} 
\begin{eqnarray}
&\sup_{x\in \mathbb{R}}&\mathbb{P}(S_{n}\in \lbrack x,x+z])\leq
C(n)z^{c}\qquad (\forall z>0)  \label{conc} \\
&&\mbox{for suitable constants $c,C(n)>0$,\quad if $n=1,2,\ldots$}  \notag
\end{eqnarray}
We also consider a bandwidth function $b\equiv b_{n}\rightarrow \infty $
such that $\lim_{n\rightarrow \infty }n/b=\infty $ and two subsampling
schemes 
\begin{eqnarray}
Y_{b,i} &=&(X_{i+1},\ldots ,X_{i+b}),\qquad N=n-b,\qquad 
\mbox{
overlapping
samples},  \label{over} \\
Y_{b,i} &=&(X_{(i-1)b+1},\ldots ,X_{ib}),\qquad N=\frac{n}{b},\qquad 
\mbox{
non-overlapping samples}.  \label{nonover2}
\end{eqnarray}
Then we introduce a smooth and a rough subsampling estimator for $\mathbb{K}$
\begin{eqnarray}
\widetilde{\mathbb{K}}_{b,n}(x)=\,\frac{1}{N}\sum_{i=0}^{N-1}\varphi \left( 
\frac{s_{b}(Y_{b,i})-x}{\epsilon _{n}}\right) , &&%
\mbox{smooth subsampled
statistics,}  \label{smooth2} \\
\widehat{\mathbb{K}}_{b,n}(x)=\,\frac{1}{N}\sum_{i=0}^{N-1}%
\mbox{I\hspace{-.6em}1}\,(s_{b}(Y_{b,i})\leq x), &&%
\mbox{rough subsampled
statistics,}  \label{rough}
\end{eqnarray}
where $\mbox{I\hspace{-.6em}1}$ is the indicator function. Here, ${\epsilon
_{n}}\downarrow 0$ and $\varphi $ is the non-increasing continuous function
such that $\varphi =1$ or 0 according to $x\leq 0$ or $x\geq 1$ and which is
affine between $0$ and $1$. From the convergent statistics assumption (\ref
{con}), one easily checks that the bias of our first estimator is bounded
the following way: 
\begin{equation}
\sup_{x\in \mathbb{R}}\Big|\mathbb{E}[\widetilde{\mathbb{K}}_{b,n}(x)]-%
\mathbb{K}(x)\Big|\leq r_{b}+\epsilon _{n}\Vert \mathbb{K}^{\prime }\Vert
_{\infty }.  \label{bias}
\end{equation}

\begin{Rem}[discussing assumptions]
\label{remhyp}\strut 

\begin{itemize}
\item  \emph{The rough subsampler (\ref{rough}) is the usual one. However in
order to derive uniform \textit{a.s.} convergence this estimator will need
the stronger concentration condition (\ref{conc}), due to the specific
problems related to weak dependence. Indeed this estimate is based on
indicators which are not Lipschitz functions and bounds for covariances are
more hard to handle. Besides the simple stationary Markov case for which
existence of a bounded transition probability density is enough to assert
that $c=1$, examples for which those concentration conditions are proved may
be found in Doukhan and Wintenberger (2007, 2008) \cite{DW} and \cite{DW2}.}

\item  \emph{The two techniques of subsampling developed here are based on
overlap\-ping or non-overlap\-ping samples; it is clear that (\ref{over}) is
much more economic in terms of the sample size $n$ since the corresponding
sum runs over $(n-b)$ indices while this number is only $n/b$ in case (\ref
{nonover2}). However the latter condition assumes less restrictive weak
dependence, since the involved $b-$samples are more distant.}
\end{itemize}
\end{Rem}

In order to prove either uniform strong or weak laws of large numbers, we
aim at bounding the $p-$th centered moments of $\widetilde{\mathbb{K}}%
_{b,n}(x)$ and $\widehat{\mathbb{K}}_{b,n}(x)$ defined as 
\begin{equation*}
\widetilde{\Delta }_{b,n}^{(p)}(x)=\left| \mathbb{E}\left[ \widetilde{%
\mathbb{K}}_{b,n}(x)-\mathbb{E[}\widetilde{\mathbb{K}}_{b,n}(x)]\right]
^{p}\right| ,\qquad \widehat{\Delta }_{b,n}^{(p)}(x)=\left| \mathbb{E}\left[ 
\widehat{\mathbb{K}}_{b,n}(x)-\mathbb{E[}\widehat{\mathbb{K}}_{b,n}(x)]%
\right] ^{p}\right| .
\end{equation*}
Borel-Cantelli lemma will then allow to conclude.

For simplicity we set the notation $L(b)={\mathrm{Lip}\,}s_{b}.\label{L(b)}$
Moreover, for two sequences $a\equiv (a_{n})$ and $b\equiv (b_{n})$, $a\prec
b$ means that there exists a positive constant $c$ such that, for all $n$, $%
a_{n}\leq bc_{n}$.

We first give results for the smooth subsampler by considering the
convergence condition (\ref{con}). Almost sure convergence is obtained to
the price of restrictive conditions that meets all the qualities required in
our framework.

\begin{theo}[Smooth subsampler]
\label{subsamp} \strut Assume the convergence assumption (\ref{con}) hold.
Let $\delta >0$ and $p\in \mathbb{N}^{\ast }$. Assume moreover that if
respectively \emph{the overlapping setting is used} and one among the
following relations hold as $n\rightarrow \infty $ 
\begin{eqnarray*}
\underline{\eta \mbox{-dependence:}} &&\sum_{t=0}^{\infty }(t+1)^{p-2}\eta
(t)<\infty ,\qquad \frac{b}{n}\left( 1\vee \frac{L(b)}{\epsilon _{n}}\right)
\prec n^{-\delta },\mbox{ or} \\
\underline{\lambda \mbox{-dependence:}} &&\sum_{t=0}^{\infty
}(t+1)^{p-2}\lambda (t)<\infty ,\qquad \frac{b}{n}\left( 1\vee \frac{L(b)}{%
\epsilon _{n}}\vee \frac{bL(b)^{2}}{\epsilon _{n}^{2}}\right) \prec
n^{-\delta },
\end{eqnarray*}
or \emph{the non-overlapping setting is used} and 
\begin{eqnarray*}
\underline{\eta \mbox{-dependence:}} &&\sum_{t=0}^{n-1}(t+1)^{p-2}\eta
(t)\prec b^{p-2},\qquad \frac{b}{n}\left( 1\vee \frac{bL(b)}{\epsilon _{n}}%
\right) \prec n^{-\delta },\mbox{ or} \\
\underline{\lambda \mbox{-dependence:}} &&\sum_{t=0}^{n-1}(t+1)^{p-2}\lambda
(t)\prec b^{p-2},\quad \frac{b}{n}\left( 1\vee \frac{bL(b)}{\epsilon _{n}}%
\vee \frac{bL(b)^{2}}{\epsilon _{n}^{2}}\right) \prec n^{-\delta }.
\end{eqnarray*}
Then 
\begin{equation*}
\widetilde{\Delta }_{b,n}^{(p)}(x)\prec n^{-[\frac{p}{2}]\delta }.
\end{equation*}
Hence, from Borel-Cantelli Lemma, if $p/2\in \mathbb{N}$ is such that $%
p\delta >2$, then 
\begin{equation*}
\sup_{x\in \mathbb{R}}\left| \widetilde{\mathbb{K}}_{b,n}(x)-\mathbb{K}%
(x)\right| \rightarrow _{n\rightarrow \infty }0\qquad \mbox{a.s.}
\end{equation*}
\end{theo}

Finally, for completeness, we give results for the rough subsampler by
considering successively the convergence condition (\ref{con}) and the
concentration condition (\ref{conc}).

\begin{theo}[Rough subsampler under condition (\ref{con})]
\label{rsubsamp2} \strut Assume that the convergence assumption (\ref{con})
holds. 
If respectively \emph{the overlapping setting is used} and one among the
following relations hold 
\begin{eqnarray*}
\underline{\eta \mbox{-dependence:}} &\displaystyle\sum_{t=0}^{\infty }\eta
(t)^{\frac{1}{2}}<\infty ,&\lim_{n\rightarrow \infty }\frac{b}{n}\left(
1\vee \frac{L({b})}{\sqrt{{b}}}\right) =0 \\
\underline{\lambda \mbox{-dependence:}} &\displaystyle\sum_{t=0}^{\infty
}\lambda (t)^{\frac{2}{3}}<\infty ,&\lim_{n\rightarrow \infty }\frac{b}{n}%
\left( 1\vee \left( \frac{L(b)^{4}}{{b}}\right) ^{\frac{1}{3}}\vee \left( 
\frac{L({b})}{{b}}\right) ^{\frac{2}{3}}\right) =0
\end{eqnarray*}
or \emph{the non-overlapping setting is used} and 
\begin{eqnarray*}
\underline{\eta \mbox{-dependence:}} &\displaystyle\sum_{t=0}^{\infty }\eta
(t)^{\frac{1}{2}}<\infty ,&\lim_{n\rightarrow \infty }\frac{{b}}{n}\left(
1\vee \sqrt{{b}}L({b})\right) =0 \\
\underline{\lambda \mbox{-dependence:}} &\displaystyle\sum_{t=0}^{\infty
}\lambda (t)^{\frac{2}{3}}<\infty ,&\lim_{n\rightarrow \infty }\frac{{b}}{n}%
\left( 1\vee (bL({b})^{2})^{\frac{2}{3}}\vee \left( {b}L({b})^{2}\right) ^{%
\frac{1}{3}}\right) =0,
\end{eqnarray*}
then $\lim_{n\rightarrow \infty }\widehat{\Delta }_{b,n}^{(2)}(x)=0$ and 
\begin{equation*}
\lim_{n\rightarrow \infty }\sup_{x\in \mathbb{R}}\left| \widehat{\mathbb{K}}%
_{b,n}(x)-\mathbb{K}(x)\right| =0,\quad \mbox{in probability.}
\end{equation*}
\end{theo}

\begin{theo}[Rough subsampler under condition (\ref{conc})]
\label{rsubsamp} \strut 

Assume that the concentration assumption (\ref{conc}) holds. Let $\delta >0$
and $p\in \mathbb{N}^{\ast }$. Assume moreover that if respectively \emph{%
the overlapping setting is used} and one among the following relations hold,
as $n\rightarrow \infty $,\bigskip 

\noindent \underline{$\eta $-dependence:}\qquad $\displaystyle%
\sum_{t=0}^{\infty }(t+1)^{p-2}\eta (t)^{\frac{2+c}{1+c}}<\infty ,$ 
\begin{equation*}
\frac{b}{n}\left( 1\vee \left( C(b)b^{-1}L(b)^{c}\right) ^{\frac{1}{1+c}%
}\vee \left( bL(b)^{2+c}\right) ^{\frac{1}{1+c}}\right) \prec n^{-\delta },
\end{equation*}
\underline{$\lambda $-dependence:}\qquad $\displaystyle\sum_{t=0}^{\infty
}(t+1)^{p-2}\lambda (t)^{\frac{1+c}{2+c}}<\infty $, 
\begin{equation*}
\frac{b}{n}\left( 1\vee (C(b)^{2}b^{c-2}L(b)^{2c})^{\frac{1}{2+c}}\vee
(C(b)b^{-2}L(b)^{c})^{\frac{1}{2+c}}\vee (C(b)b^{c-2}L(b)^{2c})^{\frac{1}{2+c%
}}\right) \prec n^{-\delta }
\end{equation*}
or \emph{the non-overlapping setting is used}\bigskip 

\noindent \underline{$\eta $-dependence:} \qquad $\displaystyle%
\sum_{t=0}^{n-1}(t+1)^{p-2}\eta (t)^{\frac{2+c}{1+c}}\prec b^{p-2},$ 
\begin{equation*}
\frac{b}{n}\left( 1\vee \left( C(b)(bL(b))^{c}\right) ^{\frac{1}{1+c}}\vee
(bL(b))^{\frac{2+c}{1+c}}\right) \prec n^{-\delta },
\end{equation*}
\underline{$\lambda $-dependence:}\qquad $\displaystyle%
\sum_{t=0}^{n-1}(t+1)^{p-2}\lambda (t)^{\frac{1+c}{2+c}}\prec b^{p-2}$ 
\begin{equation*}
\frac{b}{n}\Big(1\vee (C(b)(bL(b))^{c})^{\frac{2}{2+c}}\vee
(C(b)(bL(b))^{c})^{\frac{1}{2+c}}\vee (C(b)(bL(b))^{2c})^{\frac{1}{2+c}}\Big)%
\prec n^{-\delta }.
\end{equation*}
Then 
\begin{equation*}
\widehat{\Delta }_{b,n}^{(p)}(x)\prec n^{-[\frac{p}{2}]\delta }.
\end{equation*}
Hence, if $p/2\in \mathbb{N}$ is such that $p\delta >2$, then 
\begin{equation*}
\sup_{x\in \mathbb{R}}\left| \widehat{\mathbb{K}}_{b,n}(x)-\mathbb{K}%
(x)\right| \rightarrow _{n\rightarrow \infty }0\qquad \mbox{a.s.}
\end{equation*}
\end{theo}

The rough subsampler needs a very strong concentration assumption and
excessively intricate weak dependence conditions for uniform strong
consistency. Such conditions are definitely hard to derive in the most
general settings.

\begin{Rem}
\label{remhypth} \strut 

\begin{itemize}
\item  \emph{For Theorem \ref{subsamp} and \ref{rsubsamp}, if one of the
above mentioned relations hold with }$p=2$\emph{, then we obtain the same
result but only with respect to the convergence in probability.}

\item  Choosing the procedure. \emph{The overlapping frame yields a more
expensive procedure, in terms of the assumptions of the bandwidth. A strange
feature of the results is that, for the \emph{a.s.} convergence case where
moments with high order need to be calculated, that weak dependence
assumptions are weaker in this case. }


\item  Monitoring the smoothing parameter $\epsilon _{n}$.\emph{\ A bit more
may be found from the previous result. The square of bias of our statistics
is indeed given as $O(\epsilon _{n}^{2}+r_{b}^{2})$ while now the order of
variance is respectively 
\begin{eqnarray*}
\mbox{\underline{Overlapping}}\qquad \qquad \qquad  &\eta \mbox{-dependence:}%
&\frac{b}{n}\left( 1\vee \frac{L(b)}{\epsilon _{n}}\right) , \\
&\lambda \mbox{-dependence:}&\frac{b}{n}\left( 1\vee \frac{L(b)}{\epsilon
_{n}}\vee \frac{bL(b)^{2}}{\epsilon _{n}^{2}}\right) , \\
\mbox{\underline{Non-overlapping}}\qquad \qquad \  &\eta \mbox{-dependence:}&%
\frac{b}{n}\left( 1\vee \frac{bL(b)}{\epsilon _{n}}\right) , \\
&\lambda \mbox{-dependence:}&\frac{b}{n}\left( 1\vee \frac{bL(b)}{\epsilon
_{n}}\vee \frac{bL(b)^{2}}{\epsilon _{n}^{2}}\right) .
\end{eqnarray*}
Hence eg. under $\eta -$dependence a reasonable choice for this parameter is 
}$\epsilon _{n}$\emph{$=(bL(b)/n)^{1/3}$ (resp. $(b^{2}L(b)/n)^{1/3}$ for
the non-overlapping case). Notice however that the order of the quadratic
approximation of our subsampler is always bounded by $\max \{r_{b}^{2},b/n\}$%
.\newline
If the process is centered, the CLT writes with the statistics $%
t_{b}(x_{1},\ldots ,x_{b})=(x_{1}+\cdots +x_{b})/\sqrt{b}$; here $%
L(b)=b^{-1/2}$ and we choose respectively }$\epsilon _{n}$\emph{$%
=b^{1/6}n^{-1/3}$, or $b^{1/2}n^{-1/3}$. }

\item  Uniformity. \emph{The fact that $\varphi $ is monotonous is essential
in order to derive uniform convergence of those subsamplers; it indeed
allows to use the standard variant of Dini Theorem.}

\item  \emph{\emph{Confidence bands.} For a statistical validation of the
technique, a CLT theorem is also required. A first step for such a central
limit theorem is to precise the asymptotic variance. For simplicity we shall
mention hard subsampling of a convergent sequence of statistics (\ref{con}).
In this case $\lim_{b}\mathbb{K}_{b}(x)=\mathbb{K}(x)$ and the variance of $%
\widehat{\mathbb{K}}_{b,n}$ should involve also $t_{b,i}(x)=\mbox{Cov}\,(%
\mbox{I\hspace{-.6em}1}_{s_{b}(Y_{n,0})\leq x},\mbox{I\hspace{-.6em}1}%
_{s_{b}(Y_{n,i})\leq x})$. As in the proofs a concentration assumption leads
to $|t_{b,j}(x)|\leq \mbox{const}\cdot \eta ^{c}(j-b)$ (resp. $\leq %
\mbox{const}\cdot \lambda ^{c}(j-b)$) in the case $j>b$ and for a suitable
constant $0<c<1$; for the non-overlapping scheme we only need $j\geq 1$ and $%
j-b$ is now replaced by $(j-1)b$. The claim is now that} 
\begin{equation*}
n\mbox{Var}\,\widehat{\mathbb{K}}_{b,n}(x)\rightarrow
\sum_{k}t_{|k|}(x),\quad t_{0}(x)=\mathbb{K}(x)(1-\mathbb{K}(x)),\
t_{j}(x)=\lim_{b\rightarrow \infty }t_{b,j}(x).
\end{equation*}
\emph{An analogue result in the nonverlapping case writes more simply} 
\begin{equation*}
\frac{n}{b}\mbox{Var}\,\widehat{\mathbb{K}}_{b,n}(x)\rightarrow
t_{0}(x)+2t_{1}(x),\ t_{0}(x)=\mathbb{K}(x)(1-\mathbb{K}(x)),\
t_{1}(x)=\lim_{b\rightarrow \infty }t_{b,1}(x).
\end{equation*}
\emph{This is a first step for a CLT because in both cases the normalization
coefficient is $N$. Anyway proving a CLT involves a more precise analysis of
the situation and the use of Lindeberg method with Bernstein blocs. However,
the knowledge of this limit variance already provides a reasonable
confidence band for this estimator.}
\end{itemize}
\end{Rem}

\section{Subsampling the distribution of extremes}

Bertail \textit{et al.} (2004) \cite{BHPW} studied subsampling estimators
for distributions of diverging statistics, but imposed that the underlying
sequence is strongly mixing. We aim at generalizing their results for weakly
dependent sequences. Instead of considering the general case, we focus on
the sample maximum because we are able to give sufficient conditions such
that the normalized sample maximum converges in distribution under the weak
dependence assumption. Note however that the results can easily be
generalized provided that it is possible to compute the Lipschitz
coefficient of the function used to define the diverging statistics.

\subsection{Convergence of the sample maximum}

We first discuss conditions for convergence in distribution of the
normalized sample maximum of the weakly dependent sequence.

Let $x_{F}=\sup \{x:F(x)<1\}$ be the upper end point of $F$ and $\bar{F}%
:=1-F $. We say that the stationary distribution $F$ is in the domain
attraction of the generalized extreme value distribution with index $\gamma $%
, $-\infty <\gamma <\infty $, if there exists a positive and measurable
function $g$ such that for $1+\gamma x>0$%
\begin{equation*}
\lim_{u\rightarrow x_{F}}\bar{F}(u+xg(u))/\bar{F}(u)=(1+\gamma x)^{-1/\gamma
}\text{.}
\end{equation*}
Then there exist sequences $(u_{n})_{n\geq 1}$ and $(v_{n})_{n\geq 1}$ such
that $u_{n}>0$ and 
\begin{equation}
\lim_{n\rightarrow \infty }F^{n}(w_{n}(x))=G_{\gamma }(x):=\left\{ 
\begin{tabular}{ll}
$\exp (-(1+\gamma x)_{+}^{-1/\gamma })$ & if $\gamma \neq 0,$ \\ 
$\exp (-\exp (-x))$ & if $\gamma =0,$%
\end{tabular}
\right.  \label{convclass}
\end{equation}
where $w_{n}(x)=x/u_{n}+v_{n}$. Let $q\left( t\right) =F^{\leftarrow }\left(
1-t^{-1}\right) $ where $F^{\leftarrow }$ is the generalised inverse of $F$.
Then $(u_{n})_{n\geq 1}$ and $(v_{n})_{n\geq 1}$ can be chosen as 
\begin{eqnarray*}
v_{n} &=&q\left( n\right) \\
u_{n}^{-1} &=&\left\{ 
\begin{tabular}{ll}
$(-\gamma )(x_{F}-q\left( n\right) )$ & if $\gamma <0,$ \\ 
$q\left( ne\right) -q\left( n\right) $ & if $\gamma =0,$ \\ 
$\gamma q\left( n\right) $ & if $\gamma >0.$%
\end{tabular}
\right.
\end{eqnarray*}

Let us introduce the extremal dependence coefficient 
\begin{equation*}
\beta _{n,l}=\sup |\mathbb{P}\left( X_{i}\leq w_{n}(x),i\in A\cup B\right) -%
\mathbb{P}\left( X_{i}\leq w_{n}(x),i\in A\right) \mathbb{P}\left( X_{i}\leq
w_{n}(x),i\in B\right) |\text{,}
\end{equation*}
where the sets $A$\ and $B$\ are such that : $A\subset \left\{
1,\ldots,k\right\} $, $B\subset \left\{ k+l,\ldots,n\right\} $, and $1\leq
k\leq n-l$.

O'Brien (1987) \cite{Ob} gave sufficient conditions such that the normalized
sample maximum $u_{n}(M_{n}-v_{n})$ converges in distribution when the
stationary distribution is in the domain attraction of some extreme value
distributions.

\begin{theo}
\label{Obrienth} \cite{Ob} Assume that $F$ is in the domain attraction of
the extreme value distribution with index $\gamma $. Let $\left(
a_{n}\right) $ be a sequence of positive integers such that $a_{n}=o\left(
n\right) $ as $n\rightarrow \infty $ and 
\begin{equation}
\lim_{n\rightarrow \infty }\frac{\mathbb{P}(M_{a_{n}}>w_{n}(x))}{a_{n}\bar{F}%
(w_{n}(x))}=\theta \in (0,1].  \label{theta}
\end{equation}
Assume that there exists a sequence $\left( l_{n}\right) $ of positive
integers such that 
\begin{equation}
l_{n}=o\left( a_{n}\right) \text{\qquad and\qquad }\frac{n}{a_{n}}\beta
_{n,l_{n}}\rightarrow _{n\rightarrow \infty }0\text{\qquad as }n\rightarrow
\infty .  \label{CondLe}
\end{equation}
Then 
\begin{equation}
\lim_{n\rightarrow \infty }\mathbb{P}(M_{n}\leq w_{n}(x))=\left\{ 
\begin{tabular}{ll}
$\exp (-\theta (1+\gamma x)_{+}^{-1/\gamma })$ & if $\gamma \neq 0,$ \\ 
$\exp (-\theta \exp (-x))$ & if $\gamma =0.$%
\end{tabular}
\right.   \label{Maxtheta}
\end{equation}
\end{theo}

The constant $\theta $ is referred to as the extremal index of $X$ (see
Leadbetter \textit{et al.} (1983) \cite{LLR}). Note that any $l_{n}=o\left(
n\right) $ such that $\beta _{n,l_{n}}\rightarrow 0$ as $n\rightarrow \infty 
$ can be used in constructing a sequence $a_{n}$ such that (\ref{CondLe}) is
satisfied by taking $a_{n}$ equal to the integer part of $\max (n\beta
_{n,l_{n}}^{1/2},(nl_{n})^{1/2})$. The condition $\beta
_{n,l_{n}}\rightarrow 0$ as $n\rightarrow \infty $ is known as the $D(w_{n})$
condition (see Leadbetter (1974) \cite{lea74}). \medskip

We provide an equivalent theorem when $X$ is assumed to be either $\eta $ or 
$\lambda $-weakly dependent.

\begin{theo}
\label{deftheta}Assume that $F$ is an absolutely continuous distribution in
the domain of attraction of the extreme value distribution with index $%
\gamma $ such that for $1+\gamma x>0$ 
\begin{equation}
\lim_{n\rightarrow \infty }\frac{\partial }{\partial x}F^{n}(w_{n}(x))=\frac{%
\partial }{\partial x}G_{\gamma }(x).  \label{Condder}
\end{equation}
Let $\left( a_{n}\right) $ be a sequence of positive integers such that $%
a_{n}=o\left( n\right) $ as $n\rightarrow \infty $ and (\ref{theta}) holds.
Assume that there exists a sequence $\left( l_{n}\right) $ of positive
integers such that $l_{n}=o\left( a_{n}\right) $ ($n\rightarrow \infty $).
If $X$ is $\eta $-weakly dependent and 
\begin{equation*}
\frac{n}{a_{n}}\left( n\eta (l_{n})u_{n}\right) ^{1/2}\rightarrow
_{n\rightarrow \infty }0,
\end{equation*}
or if $X$ is $\lambda $-weakly dependent and 
\begin{equation*}
\frac{n}{a_{n}}\left( \left[ n\lambda (l_{n})u_{n}\right] ^{1/2}\vee \left[
na_{n}\lambda (l_{n})u_{n}^{2}\right] ^{1/3}\right) \rightarrow
_{n\rightarrow \infty }0,
\end{equation*}
then (\ref{Maxtheta}) holds.
\end{theo}

\subsection{Subsampling the distribution of the normalized sample maximum}

\label{subsamex}

Consider the sequence of extreme statistics 
\begin{equation*}
M_{n}=m_{n}(X_{1},\ldots ,X_{n})=\max_{1\leq i\leq n}X_{i}.
\end{equation*}
Set $\mathbb{H}_{n}(x)=\mathbb{P}(M_{n}\leq x)$. Restate the smooth
subsampling estimates for non-normalized extremes by 
\begin{equation}
\widetilde{\mathbb{H}}_{b,n}(x)=\frac{1}{N}\sum_{i=0}^{N-1}\varphi \left( 
\frac{m_{b}(Y_{b,i})-x}{\epsilon _{n}}\right) .  \label{smooth_div}
\end{equation}
Assume, under the assumption of Theorem \ref{deftheta}, that (\ref{con})
adapted to normalized extremes holds, i.e. 
\begin{equation*}
r_{n}=\sup_{x\in \mathbb{R}}|\mathbb{H}_{n}(w_{n}(x))-\mathbb{H}%
(x)|\rightarrow _{n\rightarrow \infty }0.
\end{equation*}
where $\mathbb{H=}G_{\gamma }^{\theta }$.

Following the lines of Bertail \textit{et al.} (2004) \cite{BHPW}, we have
to impose conditions on the median and the distance between two quantiles of
the limiting distribution in order to be able to identify it. The median of
the limiting distribution to estimate is assumed to be equal to $0$ and the
distance between the quantiles is assumed to be equal to $1$. Fix $%
0<t_{1}<t_{2}<1$. Then the normalizing sequences can be estimated by 
\begin{equation}
\widetilde{v}_{b,n}=\widetilde{\mathbb{H}}_{b,n}^{\leftarrow }\left( \frac{1%
}{2}\right) ,\qquad \widetilde{u}_{b,n}=\left| \widetilde{\mathbb{H}}%
_{b,n}^{\leftarrow }(t_{2})-\widetilde{\mathbb{H}}_{b,n}^{\leftarrow
}(t_{1})\right| ^{-1}.  \label{normcoeff1}
\end{equation}
Let $C=\mathbb{H}^{\leftarrow }(t_{2})-\mathbb{H}^{\leftarrow }(t_{1})$.
Using that here $L(b)={\mathrm{Lip}\,}m_{b}=1$, we derive from Theorem \ref
{subsamp} and Theorem 4 in Bertail \textit{et al.} (2004) \cite{BHPW} the
following theorem.

\begin{theo}
\label{subsampEX} Assume that the conditions of Theorem \ref{deftheta} hold.
Let $\delta >0$ and $p\in \mathbb{N}^{\ast }$.

The relation $\left| \mathbb{E}\left[ \widetilde{\mathbb{H}}_{b,n}(w_{b}(x))-%
\mathbb{E[}\widetilde{\mathbb{H}}_{b,n}(w_{b}(x))]\right] ^{p}\right| \prec
n^{-[\frac{p}{2}]\delta }$ holds if we assume that $\lim_{n\rightarrow
\infty }\epsilon _{n}u_{b}=0$ and respectively (as $n\rightarrow \infty $)
that:

\begin{itemize}
\item  \emph{in the overlapping case}, \newline
\underline{$\eta $-weak dependence,}\qquad $\displaystyle\sum_{t=0}^{\infty
}(t+1)^{p-2}\eta (t)<\infty ,\quad \frac{b}{n^{1-\delta }\epsilon _{n}}\prec
1,$ or \newline
\underline{$\lambda $-weak dependence,}\qquad $\displaystyle%
\sum_{t=0}^{\infty }(t+1)^{p-2}\lambda (t)<\infty ,\quad \frac{b}{%
n^{1-\delta }\epsilon _{n}^{2}}\prec 1,$

\item  \emph{in the non-overlapping case}, \newline
\underline{$\eta $-weak dependence,}\qquad $\displaystyle%
\sum_{t=0}^{n-1}(t+1)^{p-2}\eta (t)\prec b^{p-2},\quad \frac{b^{2}}{%
n^{1-\delta }\epsilon _{n}}\prec 1,$ or \newline
\underline{$\lambda $-weak dependence,}\qquad $\displaystyle%
\sum_{t=0}^{n-1}(t+1)^{p-2}\lambda (t)\prec b^{p-2},\quad \frac{b^{2}}{%
n^{1-\delta }\epsilon _{n}^{2}}\prec 1.$
\end{itemize}

Hence, if $p/2\in \mathbb{N}$ is such that $p\delta >2$, then 
\begin{equation*}
\sup_{x\in \mathbb{R}}\left| \widetilde{\mathbb{H}}_{b,n}\left( \widetilde{v}%
_{b,n}+\frac{x}{\widetilde{u}_{b,n}}\right) -\mathbb{H}\left( \mathbb{H}%
^{\leftarrow }\left( \frac{1}{2}\right) +Cx\right) \right| \rightarrow
_{n\rightarrow \infty }0\qquad \mbox{a.s.}
\end{equation*}
\newline
\end{theo}


\section{Simulation study}

\label{SEx}

The finite sample properties of our subsampling estimators are now compared
in a simulation study. We consider both rough and smooth subsampling
estimators when they are computed with the overlapping or non overlapping
schemes.

Sequences of length $n=2$,$000$ and $n=5,000$ have been simulated from the
first-order autoregressive process of Example (\ref{FOS}) 
\begin{equation*}
X_{t}=\frac{1}{r}(X_{t-1}+\varepsilon _{t}),
\end{equation*}
where $(\varepsilon _{t})_{t\in \mathbb{Z}}$ are iid and uniformly
distributed on the set $\{0,1,\ldots ,r-1\}$ and $r$ is equal to $3$. It is
well-known that the asymptotic condition (\ref{Maxtheta}) holds with $\gamma
=-1$, $\theta =r^{-1}(r-1)$, $u_{n}=n$ and $v_{n}=1-n^{-1}$. Following the
approach presented in the previous Subsection \ref{subsamex}, we have to fix
conditions on the median and two quantiles of the limiting distribution. We
choose $t_{1}=1/4$ and $t_{2}=3/4$. The limiting distribution becomes 
\begin{equation*}
\mathbb{K}(x)=e^{-\theta (1-(x-d)/c)},\qquad x\leq c+d,
\end{equation*}
where $c=\theta /\ln 3$ and $d=(\ln 2-\theta )/\ln 3$. The normalization
coefficients $\bar{u}_{n}$ and $\bar{v}_{n}$ such that 
\begin{equation*}
\lim_{n\rightarrow \infty }\mathbb{P}(\bar{u}_{n}(M_{n}-\bar{v}_{n})\leq x)=%
\mathbb{K}(x)
\end{equation*}
are given by 
\begin{equation}
\bar{v}_{n}=v_{n}-c^{-1}u_{n}^{-1}d,\qquad \bar{u}_{n}=cu_{n}.
\label{normcoff2}
\end{equation}

We first simulate a sequence of length $n=2,000$ and plot the estimators of
the limiting distribution in Figure 1. As expected, smoothing estimators
yield smoother curves. The differences between the estimators are small but
the smoothed versions need less strong assumptions for their trajectorial
convergence. 
\begin{figure}[h]
\label{graph0} \hspace{4cm} \includegraphics[scale=0.5]{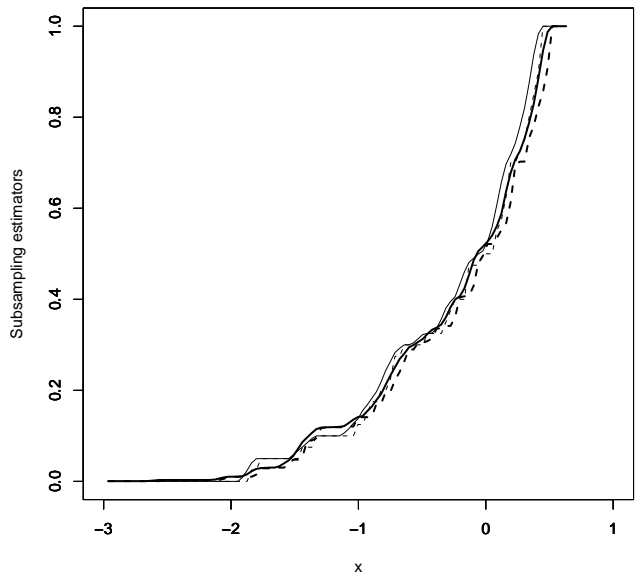}
\caption{AR(1) process - \textit{The rough (dashed) and smooth (solid)
subsampling estimators computed with the non-overlapping scheme (thin) and
with the overlapping scheme (thick) for a sequence of length }$n=2,000$%
\textit{. }$b=50$, $\protect\epsilon =0.05$.}
\end{figure}

Monte Carlo approximations to the quantiles and the means of the estimators
have been then computed from $1,000$ simulated sequences.

The properties of our rough and smooth subsampling estimators computed with
the non-overlapping scheme are shown in the two upper graphs in Figure 2.
There are very few differences between both estimators according to their
quantiles and their means. Their biases are negligible for all the value of $%
x$. The confidence intervals with level $90\%$ (gray zone) vanish when $x$
goes to $0$ because $0$ is the median of the empirical distribution, but
also the median of the asymptotic distribution. We may compare the quantiles
and the means of our estimators with those obtained when the normalization
coefficients given by (\ref{normcoeff1}) are replaced by the theoretical
normalization coefficients given by (\ref{normcoff2}) (see the two lower
graphs in Figure 2). First note that the bias become negative when $x$ is
smaller than the median. Second the confidence intervals are obviously not
equal to zero for the median but they are more narrow than the confidence
intervals of our estimators when $x$ is close to the extremal point of the
asymptotic distribution, $c+d$.

The properties of our rough and smooth subsampling estimators computed with
the overlapping scheme are shown in Figure 3. We chose the same value for $b$
as in the non-overlapping scheme and consequently the number of components
in the definition the estimators is quite larger than in the other scheme.
It follows that the empirical distribution functions given by the estimators
computed with the overlapping scheme are smoother than those of the
estimators computed with the non-overlapping scheme. The intervals
confidence are also a little bit more narrow.

Moreover note that qualitatively similar results were found when the
simulations were repeated with $n=5,000$, $b=100$, $\epsilon =0.05$.

\begin{center}
\begin{figure}[h]
\label{graph1} \hspace{3cm} %
\includegraphics[scale=0.6]{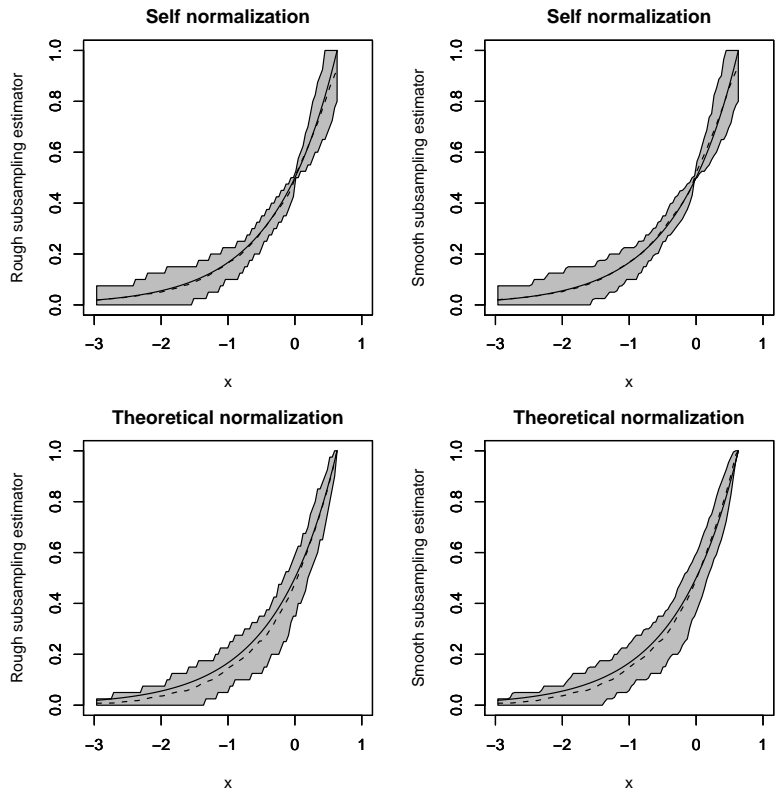}
\caption{AR(1) process - \textit{Monte Carlo approximations to the quantiles
($q_{0.05}$ and $q_{0.95}$) (gray zone) and means (dashed line) of the rough
(left) and smooth (right) subsampling estimators computed with the
non-overlapping scheme when the normalization coefficients are given by (\ref
{normcoeff1}) (top) or by (\ref{normcoff2}) (bottom). The asymptotic
distribution function, }$\mathbb{K}$\textit{, is given by the solid line. }$%
n=2,000$, $b=50$, $\protect\epsilon =0.05$.}
\end{figure}
\end{center}

\strut

Finally sequences of length $n=2$,$000$ have also been simulated from the
LARCH model with Rademacher iid inputs (\ref{Larch}) and with inputs that
have a parabolic density probability function given by $x\mapsto 0.5(1+\rho
)|x|^{\rho }$ for $x\in \left[ -1,1\right] $. Note that the Rademacher
distribution can be seen as the limit of the parabolic distribution as $\rho 
$ goes to infinity. We choose $a=0.4$. Hence the process is weak dependent
but not strong mixing when the inputs have a Rademacher distribution, and it
is strong mixing when the distribution of the inputs is absolutely
continuous. Neither the stationary distribution, nor the extremal behavior
of the processes are known. Note however that the end points of the
stationary distributions are finite.

We perform simulations and use our estimators. Results are given in Figure
4. The shapes of the empirical distribution functions given by the
estimators are different for the two processes (in particular for the large
values of $x$). As far as we can see, the generalized extreme value
distribution with a negative index could be a good choice to model the
distribution of the maximum of the process with absolutely continuous inputs
but not to model the distribution of the maximum of the process with
Rademacher inputs. The study of the extremal behavior of these processes are
intricate and left for future work.

\begin{center}
\begin{figure}[h]
\label{graph2} \hspace{3cm} %
\includegraphics[scale=0.6]{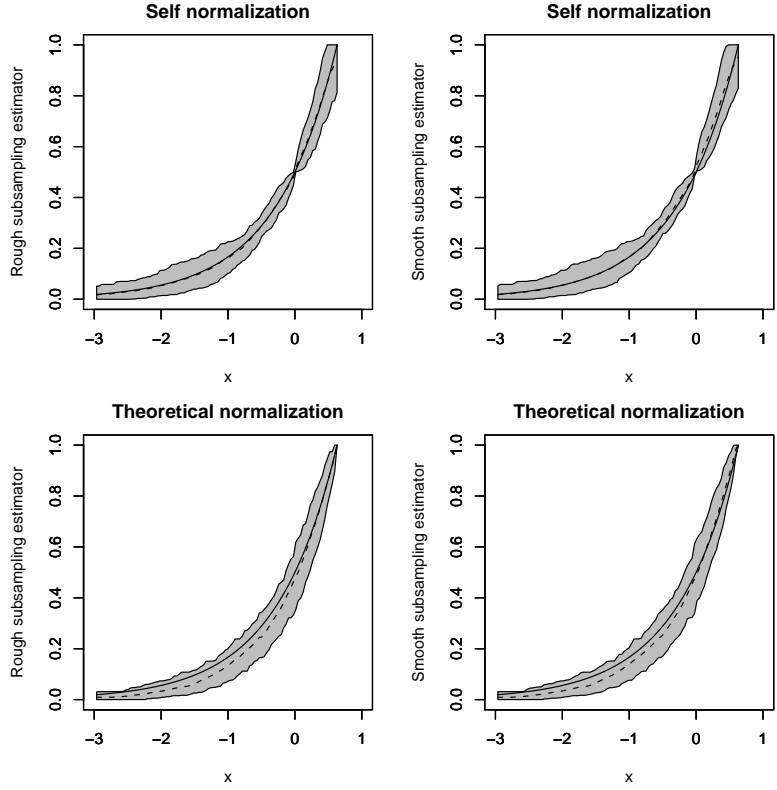}
\caption{AR(1) process - \textit{Monte Carlo approximations to the quantiles
($q_{0.05}$ and $q_{0.95}$) (gray zone) and means (dashed line) of the rough
(left) and smooth (right) subsampling estimators computed with the
overlapping scheme when the normalization coefficients are given by (\ref
{normcoeff1}) (top) or by (\ref{normcoff2}) (bottom). The asymptotic
distribution function, }$\mathbb{K}$\textit{, is given by the solid line. }$%
n=2,000$, $b=50$, $\protect\epsilon =0.05$.}
\end{figure}

\begin{figure}[h]
\label{graph3} \includegraphics[scale=0.4]{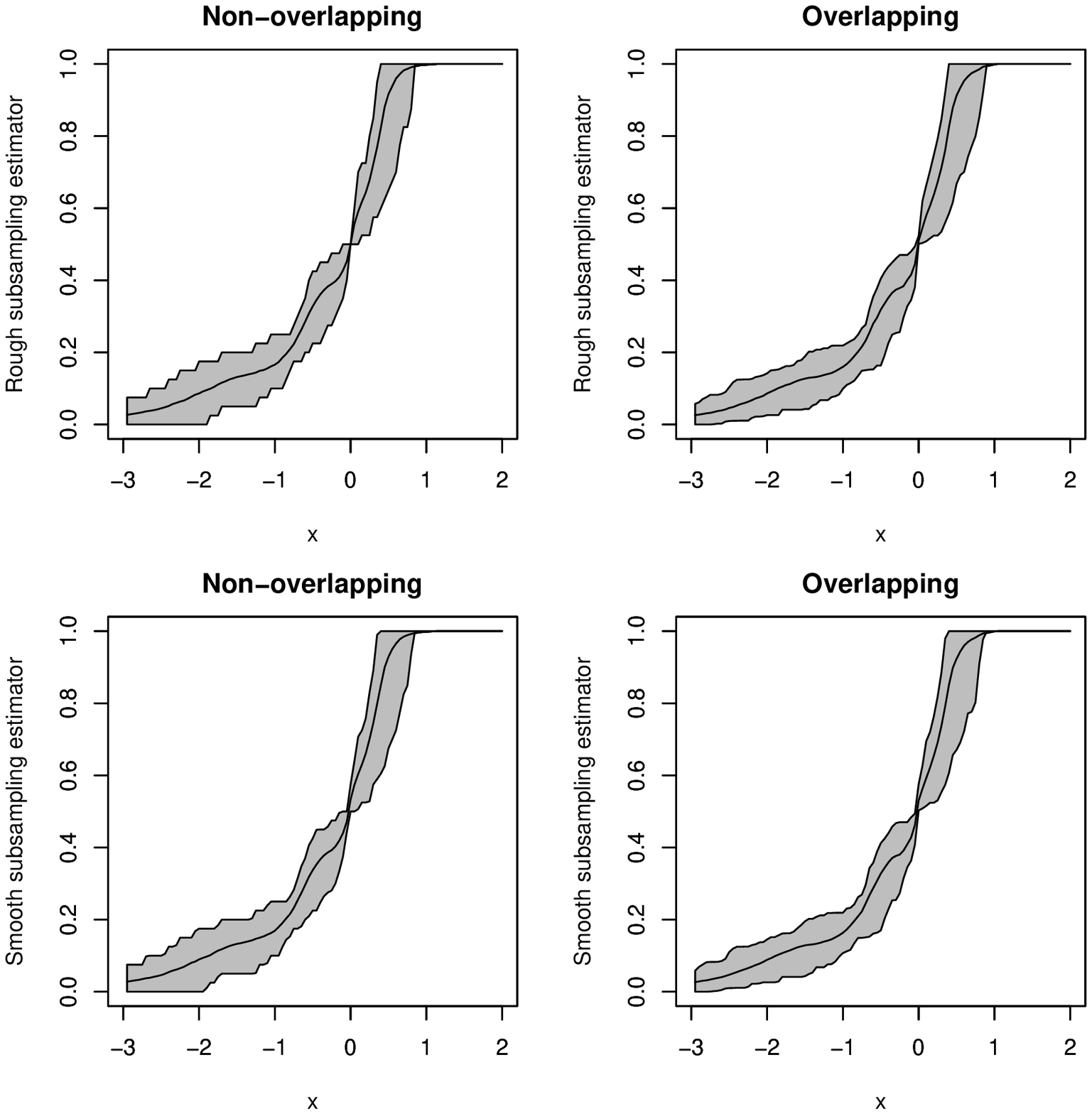}\hspace{3cm}%
\includegraphics[scale=0.4]{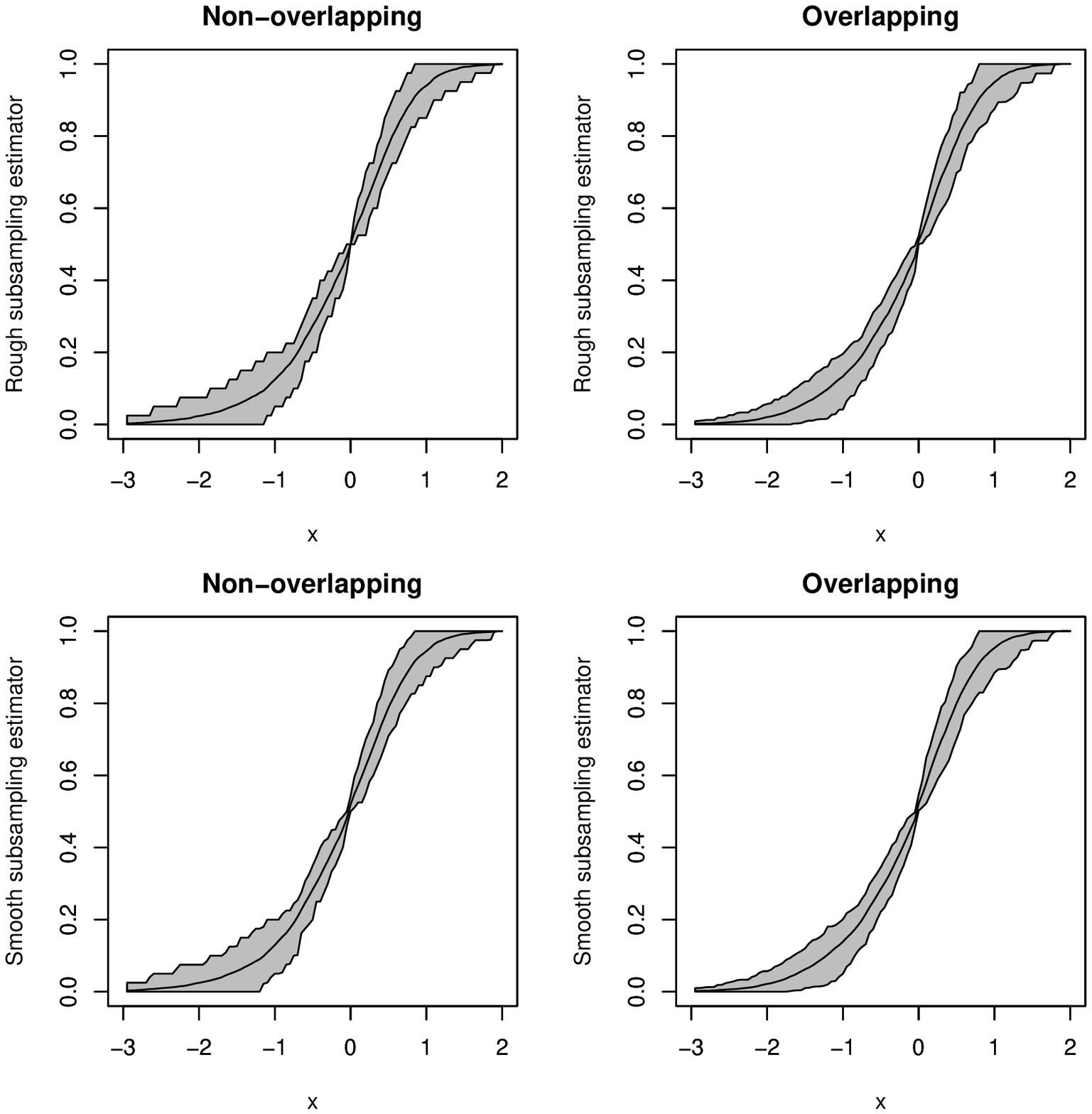}
\caption{LARCH processes with Rademacher inputs (the four graphs at left)
and with parabolic inputs ($\protect\rho =4$) (the four graphs at right) - 
\textit{Monte Carlo approximations to the quantiles ($q_{0.05}$ and $q_{0.95}
$) (gray zone) and means (solid line) of the rough (top) and smooth (bottom)
subsampling estimators computed with the non-overlapping scheme (left) and
with the overlapping scheme (right). }$n=2,000$, $b=50$, $\protect\epsilon %
=0.05$.}
\end{figure}
\end{center}

\section{Proofs}

\label{proofs}

\subsection{Proofs for smooth subsampling}

A bound of the expression $\widetilde\Delta_{b,n}^{(p)}(x)$ is closely
related to the coefficients defined for $1\le q\le p$ as: 
\begin{equation*}
C_{b,q}(r)=\sup\left|\mbox{Cov}\, \left(Z_{i_1}\cdots
Z_{i_k},Z_{i_{k+1}}\cdots Z_{i_q}\right)\right|
\end{equation*}
where the supremum refers to such indices with $1\le k<q$, $i_1\le \cdots\le
i_q$ satisfy $i_{k+1}-i_k=r$ and $Z_i=\varphi\left(\frac{s_{b}(Y_{b,i})-x}{%
\epsilon_n} \right)-\mathbb{E}\varphi\left(\frac{s_{b}(Y_{b,i})-x}{\epsilon_n%
}\right)$ is a centered rv. Then setting 
\begin{equation*}
A_{b,q}(N)=\frac1{N^q}\sum_{1\le i_1 \le \cdots\le i_q\le N}\left|\mathbb{E}
Z_{i_1}\cdots Z_{i_q}\right|, \qquad 2\le q\le p
\end{equation*}
Doukhan and Louhichi (1999) prove that $\displaystyle
\widetilde\Delta_{b,n}^{(p)}(x)\le p!A_{b,p}(N)$, moreover: 
\begin{eqnarray*}
A_{b,p}(N)&\le& B_{b,p}(N)+\sum_{q=2}^{p-2}A_{b,q}(N)A_{b,p-q}(N) \\
B_{b,q}(N)&=&\frac{q-1}{N^{q-1}}\sum_{r=0}^{N-1}(r+1)^{q-2}C_{b,q}(r),
\qquad 2\le q\le p
\end{eqnarray*}

\begin{Lemma}
\label{gen_b} Let $p,q,b,N$ be integers and $\beta (b,N)\leq 1$, we assume
that for all $2\leq q\leq p$ there exists a constant $c_{q}\geq 0$ such that 
$B_{b,q}(N)\leq c_{q}\beta ^{\frac{q}{2}}(b,N)$. Then there exists a
constant $C_{p}\geq 0$ only depending on $p$ and $c_{1},\dots ,c_{p}$ such
that $A_{b,p}(N)\leq C_{p}\beta ^{\lbrack \frac{p}{2}]}(b,N)$.
\end{Lemma}

\begin{proofl}
The result is the assumption if $p=2$ because $A_{b,2}(N)\leq B_{b,2}(N)$.
If now the result has been proved for each $q<p$ the relation $[\frac{p}{2}%
]\leq \lbrack \frac{q}{2}]+[\frac{p-q}{2}]$ completes the proof because $%
\beta (b,N)\leq 1$.
\end{proofl}

A covariance $\mbox{Cov}\,\big(f(Y_{i_1},\ldots, Y_{i_u}),g(Y_{j_1},\ldots,
Y_{j_v})\big)$ writes respectively as 
\begin{equation*}
\left\{ 
\begin{array}{lr}
\mbox{Cov}\,\Big(f_b\left(\left(X_{i_h+k}\right)_{{\tiny \left\{ 
\begin{array}{l}
1\le h \le u \\ 
1\le k\le b
\end{array}
\right.}}\right), & g_b\left(\left(X_{j_{h^{\prime }}+k^{\prime }}\right)_{%
{\tiny \left\{ 
\begin{array}{l}
1\le h^{\prime }\le u \\ 
1\le k^{\prime }\le b
\end{array}
\right.}}\right)\Big) \\ 
\mbox{Cov}\,\Big(f_b\left( \left(X_{(i_h-1)b+k}\right)_{{\tiny \left\{ 
\begin{array}{l}
1\le h \le u \\ 
1\le k\le b
\end{array}
\right.}} \right), & g_b\left(\left(X_{(j_{h^{\prime }}-1)b+k^{\prime
}}\right)_{{\tiny \left\{ 
\begin{array}{l}
1\le h^{\prime }\le u \\ 
1\le k^{\prime }\le b
\end{array}
\right.}}\right)\Big)
\end{array}
\right.
\end{equation*}
for suitable functions $f_b,g_b$ depending if the considered setting is the
overlapping one or not. Moreover ${\mathrm{Lip } \, } f_b\le {\mathrm{Lip }
\, } f$ which proves that if the dependence coefficients relative to the
sequences ${\ X}=(X_i)_{i\in \mathbb{Z}}$ are denoted by $\eta_{
X}(r)=\eta(r)$ and $\eta_{Y_b}(r)$, then we get the elementary lemma

\begin{Lemma}[Heredity]
\label{her_b} Assume that the stationary sequence ${X}=(X_{i})_{i\in \mathbb{%
Z}}$ is weakly dependent then the same occurs for ${Y_{b}}=(Y_{i})_{i\in 
\mathbb{Z}}$ and:

\begin{itemize}
\item  $\eta _{Y_{b}}(r)\leq b\eta (r-b)$ if $r\geq b$ in the overlapping
case,

\item  $\lambda _{Y_{b}}(r)\leq b^{2}\lambda (r-b)$ if $r\geq b$ in the
overlapping case,

\item  $\eta _{Y_{b}}(r)\leq b\eta ((r-1)b)$ if $r\geq 1$ in the
non-overlapping case,

\item  $\lambda _{Y_{b}}(r)\leq b^{2}\lambda ((r-1)b)$ if $r\geq 1$ in the
non-overlapping case.
\end{itemize}
\end{Lemma}

\noindent In our setting we use the function $\displaystyle %
f(y_1,\ldots,y_b)=\varphi\left(\frac{s_b(y_1,\ldots,y_b)-x}{\epsilon_n}
\right)$, the covariance inequalities write here as:

\begin{Lemma}
\label{coefC_ptilde} Using conditions (\ref{over},\ref{nonover2}) and under
the respective weak dependence assumptions $\eta $ and $\lambda $ we
respectively get

\begin{itemize}
\item  in the overlapping case, $C_{b,q}(r)\prec 1$ for $r<b$ and else,
resp. 
\begin{equation*}
C_{b,q}(r)\prec \frac{bL(b)}{\epsilon _{n}}\eta (r-b),\mbox{ or }%
C_{b,q}(r)\prec \frac{bL(b)}{\epsilon _{n}}\left( 1\vee \frac{bL(b)}{%
\epsilon _{n}}\right) \lambda (r-b)
\end{equation*}

\item  in the non-overlapping case, $C_{b,q}(r)\prec 1$ for $r=0$ and else,
resp. 
\begin{equation*}
C_{b,q}(r)\prec \frac{bL(b)}{\epsilon _{n}}\eta ((r-1)b),\mbox{ or }%
C_{b,q}(r)\prec \frac{bL(b)}{\epsilon _{n}}\left( 1\vee \frac{bL(b)}{%
\epsilon _{n}}\right) \lambda ((r-1)b)
\end{equation*}
\end{itemize}
\end{Lemma}

This lemma entails the bounds:

\begin{itemize}
\item  \underline{Overlapping and $\eta $-dependent case.} We obtain
\end{itemize}

\begin{eqnarray*}
B_{b,q}(N)&\prec& \frac{1}{N^{q-1}}\sum_{r=0}^{b-1}(r+1)^{q-2}+%
\frac1{N^{q-1}} \frac{bL(b)}{\epsilon _{n}} \sum_{r=b}^{N-1}\frac{\eta (r-b) 
}{(r+1)^{2-q}} \\
&\prec& \left( \frac{b}{N}\right) ^{q-1}\left(1+ \frac{L(b)}{\epsilon _{n}}%
\sum_{t=0}^{N-b-1}\eta (t)\right) +\frac1{N^{q-1}}\frac{bL(b)}{\epsilon _{n}}%
\sum_{t=0}^{N-b-1}\frac{\eta (t)}{(t+1)^{2-q}}
\end{eqnarray*}
where the second inequality follows from the change in variable $r=t+b$. We
use here $N=n$. Now if we assume $b\prec n^{1-\delta }$ we deduce that $%
(b/n)^{q-1}\prec n^{-\frac{q}{2}\delta}$ and if $bL(b)\prec
n^{1-\delta}\epsilon _{n}$ we analogously derive that $(b/N)^{q-1}(L(b)/%
\epsilon _{n})\prec n^{-\frac{q}{2}\delta}$. Assume now that 
\begin{equation*}
\frac{b}{n}\left( 1\vee \frac{L(b)}{\epsilon _{n}}\right) \prec n^{-\delta }
\end{equation*}
this implies with $\displaystyle
\sum_{t=0}^{\infty }\left( t+1\right) ^{q-2}\eta (t)<\infty $ that $%
B_{b,q}(N)\prec n^{-\frac{q}{2}\delta}$.

\begin{itemize}
\item  \underline{Overlapping and $\lambda $-dependent case.} We obtain
\end{itemize}

\begin{eqnarray*}
B_{b,q}(N) &\prec &\frac{1}{N^{q-1}}\sum_{r=0}^{b-1}(r+1)^{q-2}+%
\frac1{N^{q-1}} \frac{bL(b)}{\epsilon _{n}}\sum_{r=b}^{N-1}\frac{\lambda
(r-b)}{(r+1)^{2-q}} \\
&&+\,\frac1{N^{q-1}}\frac{(bL(b))^{2}}{\epsilon _{n}^{2}}\sum_{r=b}^{N-1}%
\frac{\lambda (r-b) }{(r+1)^{2-q}} \\
&\prec &\left( \frac{b}{N}\right) ^{q-1} \Bigl(1+\frac{L(b)}{\epsilon _{n}}%
\sum_{t=0}^{N-b-1}\lambda (t) + \frac{bL(b)^{2}}{\epsilon _{n}^{2}}%
\sum_{t=0}^{N-b-1}\lambda (t) \Bigr) \\
&&+\,\frac1{N^{q-1}}\left(\frac{bL(b)}{\epsilon _{n}}\sum_{t=0}^{N-b-1} 
\frac{ \lambda(t)}{(t+1)^{2-q}}+\frac{(bL(b))^{2}}{\epsilon _{n}^{2}}%
\sum_{t=0}^{N-b-1}\frac{\lambda (t)}{(t+1)^{2-q}}\right)
\end{eqnarray*}
where the second inequality follows from the change in variable $r=t+b$. We
use here $N=n$. Now if we assume $b\prec n^{1-\delta }$ we deduce that $%
(b/n)^{q-1}\prec n^{-\frac{q}{2}\delta}$. Now if $(bL(b)\vee
(bL(b))^{2})\prec n^{1-\delta }(\epsilon _{n}\vee \epsilon _{n}^{2})$ we
analogously derive that $(b/N)^{q-1}(L(b)/\epsilon _{n}\vee
bL(b)^{2}/\epsilon _{n}^{2})\prec n^{-\frac{q}{2}\delta}$. Assume now that 
\begin{equation*}
\frac{b}{n}\left( 1\vee \frac{L(b)}{\epsilon _{n}}\vee \frac{bL(b)^{2}}{%
\epsilon _{n}^{2}}\right) \prec n^{-\delta }
\end{equation*}
this implies with $\displaystyle
\sum_{t=0}^{\infty }\left( t+1\right) ^{q-2}\lambda (t)<\infty $ that $%
B_{b,q}(N)\prec n^{-\frac{q}{2}\delta}$.

\begin{itemize}
\item  \underline{Non-overlapping and $\eta $-dependent case.} We obtain
\end{itemize}

\begin{eqnarray*}
B_{b,q}(N) &\prec &\frac{1}{N^{q-1}}+\frac{1}{N^{q-1}}\frac{bL(b)}{\epsilon
_{n}}\sum_{r=1}^{N-1}\frac{\eta ((r-1)b) }{(r+1)^{2-q}} \\
&\prec & \frac{1}{N^{q-1}}\Bigl(1+\frac{bL(b)}{\epsilon _{n}}%
\sum_{k=1}^{n-1}\eta (k) \Bigr)+\frac{1}{N^{q-1}}\frac{b^{3-q}L(b)}{\epsilon
_{n}}\sum_{k=1}^{n-1}\frac{\eta (k)}{k^{2-q}} \\
\end{eqnarray*}
where the second inequality follows from replacement of $k=b(r-1)$. We use
here $N=n/b$. Now if we assume $b\prec n^{1-\delta }$ we deduce that $%
(b/n)^{q-1}\prec n^{-\frac{q}{2}\delta}$ and if $b^{2}L(b)\prec n^{1-\delta
}\epsilon _{n}$ we analogously derive that $(1/N^{q-1})(bL(b)/\epsilon
_{n})\prec n^{-\frac{q}{2}\delta}$. Assume now that 
\begin{equation*}
\frac{b}{n}\left( 1\vee \frac{bL(b)}{\epsilon _{n}}\right) \prec n^{-\delta }
\end{equation*}
this implies with $\displaystyle\sum_{t=0}^{n-1}\left( t+1\right) ^{q-2}\eta
(t)<b^{q-2}$ that $B_{b,q}(N)\prec n^{-\frac{q}{2}\delta}$.

\begin{itemize}
\item  \underline{Non-overlapping and $\lambda $-dependent case.} We obtain
\end{itemize}

\begin{eqnarray*}
B_{b,q}(N) &\prec &\frac{1}{N^{q-1}}+\frac{1}{N^{q-1}}\frac{bL(b)}{\epsilon
_{n}}\sum_{r=b}^{N-1}\frac{\lambda ((r-1)b)}{(r+1)^{2-q}} +\frac{1}{N^{q-1}}%
\frac{(bL(b))^{2}}{\epsilon _{n}^{2}}\sum_{r=b}^{N-1}\frac{\lambda ((r-1)b)}{%
(r+1)^{2-q}} \\
&\prec &\left( \frac{1}{N}\right) ^{q-1}\Bigl(1+\frac{bL(b)}{\epsilon _{n}}%
\sum_{k=b}^{n-1}\lambda (k)+\frac{bL(b)^{2}}{\epsilon _{n}^{2}}%
\sum_{k=b}^{n-1}\lambda (k)\Bigr) \\
&&+\,\frac{1}{N^{q-1}}\left(\frac{b^{3-q}L(b)}{\epsilon _{n}}%
\sum_{k=b}^{n-1} \frac{\lambda (k) }{k^{2-q}}+\frac{b^{5-q}L(b)^{2}}{%
\epsilon _{n}^{2}}\sum_{k=b}^{n-1} \frac{\lambda (k)}{k ^{2-q}}\right),
\end{eqnarray*}
where the second inequality follows from replacement of $k=b(r-1)$. We use
here $N=n/b$. Now if we assume $b\prec n^{1-\delta }$ we deduce that $%
(b/n)^{q-1}\prec n^{-\frac{q}{2}\delta}$. Now if $(b^{2}L(b)\vee
(bL(b))^{2})\prec n^{1-\delta }(\epsilon _{n}\vee \epsilon _{n}^{2})$ we
analogously derive that $(1/N^{q-1})(bL(b)/\epsilon _{n}\vee
bL(b)^{2}/\epsilon _{n}^{2})\prec n^{-\frac{q}{2}\delta}$. Assume now that 
\begin{equation*}
\frac{b}{n}\left( 1\vee \frac{bL(b)}{\epsilon _{n}}\vee \frac{bL(b)^{2}}{%
\epsilon _{n}^{2}}\right) \prec n^{-\delta }
\end{equation*}
this implies with $\displaystyle \sum_{t=0}^{n-1}\left( t+1\right)
^{q-2}\lambda (t)<\left( b^{q-2} \vee b^{q-4} \right)$ that $B_{b,q}(N)\prec
n^{-\frac{q}{2}\delta}$.

\begin{Lemma}
The relation $B_{b,q}(N)\prec n^{-\frac{q}{2}\delta}$ holds in the following
cases

\begin{itemize}
\item  In the overlapping case, if we have respectively 
\begin{eqnarray*}
&\sum_{t=0}^{\infty }(t+1)^{q-2}\eta (t)<\infty ,\mbox{ and
}&\frac{b}{n}\left( 1\vee \frac{L(b)}{\epsilon _{n}}\right) \prec n^{-\delta
},\qquad \qquad \quad \qquad \qquad \quad \\
&\sum_{t=0}^{\infty }(t+1)^{q-2}\lambda (t)<\infty ,\mbox{ and }&\frac{b}{n}%
\left( 1\vee \frac{L(b)}{\epsilon _{n}}\vee \frac{bL(b)^{2}}{\epsilon
_{n}^{2}}\right) \prec n^{-\delta }.
\end{eqnarray*}

\item  In the non-overlapping case, if we have respectively 
\begin{eqnarray*}
&\sum_{t=0}^{n-1}(t+1)^{q-2}\eta (t)<b^{q-2},\mbox{ and
}&\frac{b}{n}\left( 1\vee \frac{bL(b)}{\epsilon _{n}}\right) \prec
n^{-\delta },\qquad \qquad \quad \qquad \qquad \\
&\sum_{t=0}^{n-1}(t+1)^{q-2}\lambda (t)<\left( b^{q-2} \vee b^{q-4} \right),%
\mbox{and}&\frac{b}{n}\left( 1\vee \frac{bL(b)}{\epsilon _{n}}\vee \frac{%
bL(b)^{2}}{\epsilon _{n}^{2}}\right) \prec n^{-\delta }.
\end{eqnarray*}
\end{itemize}
\end{Lemma}

This lemma together with lemma \ref{gen_b} yields the main theorem.

\subsection{Proofs for rough subsampling}

In this section we shall replace $\epsilon_n$ by some $z>0$ to be settled
later and we set $\varphi_z(t)=\varphi\left(\frac{t-x}z\right)$. We now set $%
Z_i=\mbox{I\hspace{-.6em}1}_{\{s_{b}(Y_{b,i})\le x\}}-\mathbb{P}%
(s_{b}(Y_{b,i})\le x)$ and $W_i=\varphi_{z}(s_{b}(Y_{b,i}))-\mathbb{E}%
\varphi_{z}(s_{b}(Y_{b,i} ))$. An usual trick yields: 
\begin{eqnarray*}
\left|\mbox{Cov}\, \left(Z_{i_1}\cdots Z_{i_k},Z_{i_{k+1}}\cdots
Z_{i_q}\right)\right|&\le& \left|\mbox{Cov}\, \left(W_{i_1}\cdots
W_{i_k},W_{i_{k+1}}\cdots W_{i_q}\right)\right| \\
&+&2\sum_{h=1}^p\mathbb{E}|W_{i_h}-Z_{i_h}| \,=\, U+V
\end{eqnarray*}
with $U=\left|\mbox{Cov}\, \left(W_{i_1}\cdots W_{i_k},W_{i_{k+1}}\cdots
W_{i_q}\right)\right|$ and $V= 2p\,\mathbb{P}(s_b(Y_{b,i})\in[x,x+z])$. 
\newline
A bound for $V$ does not depend on the overlapping or not overlapping case
and we get 
\begin{equation*}
\left\{ 
\begin{array}{llllll}
V & \le & 2p(G{\mathrm{Lip } \, } s_bz+r_b) & \prec & L(b)z+r_b, & 
\mbox{under
assumption } (\ref{con}) \\ 
V & \le & 2pC(b)z^c & \prec & C(b) z^c, & \mbox{ under
assumption } (\ref{conc})
\end{array}
\right.
\end{equation*}
Set here $A_{p,b,\epsilon}={bp{\mathrm{Lip } \, } s_b}/z\prec bL(b)/z$. The
bound of $U$ needs 4 cases (considered in Lemma \ref{coefC_ptilde}) with 

\begin{itemize}
\item  in the overlapping case, $\displaystyle U\prec \left\{ 
\begin{array}{ll}
\frac{bL(b)}{z}\eta (r-b), & \mbox{for }r\geq b \\ 
1, & \mbox{for }r<b
\end{array}
\right. $

\item  in the overlapping case, $\displaystyle U\prec \left\{ 
\begin{array}{ll}
\frac{bL(b)}{z}(1\vee \frac{bL(b)}{z})\lambda (r-b), & \mbox{for }r\geq b \\ 
1, & \mbox{for
}r<b
\end{array}
\right. $

\item  in the non-overlapping case, $\displaystyle U\prec \left\{ 
\begin{array}{ll}
\frac{bL(b)}{z}\eta ((r-1)b), & \mbox{for }r\geq 1 \\ 
1, & \mbox{for }r=0
\end{array}
\right. $

\item  in the non-overlapping case, $\displaystyle
U\prec \left\{ 
\begin{array}{ll}
\frac{bL(b)}{z}(1\vee \frac{bL(b)}{z})\lambda ((r-1)b), & \mbox{for }r\geq 1
\\ 
1, & \mbox{for }r=0
\end{array}
\right. $
\end{itemize}

We first derive the inequality $(t+1+b)^{q-2}\le 2^{(q-3)}\vee1\left\{
(t+1)^{q-2}+b^{q-2}\right\}$ from convexity if $q>3$ and sublinearity else,
thus: 
\begin{equation*}
(t+1+b)^{q-2}\prec (t+1)^{q-2}+b^{q-2}.
\end{equation*}
Coefficients $C_{b,q}(r)\prec \sup\{U+V\}$ may thus be bounded in all the
considered cases.

For simplicity we classify the cases with couples of numbers indicating the
fact overlapping (\ref{over}) or not (\ref{nonover2}) setting is used and
from the fact the convergence (\ref{con}) or concentration (\ref{conc}) is
assumed, which makes 4 different cases to consider). Consider the cases
under assumption (\ref{con}).

\textbf{- }\underline{$\eta$ (\ref{over},\ref{con}) case.} Note that 
\begin{eqnarray*}
C_{b,q}(r)&\prec& L(b)\Bigl(b\eta(r-b)/z+z\Bigr)+r_b\prec L(b)\sqrt{%
b\eta(r-b)}+r_b \qquad \qquad \qquad
\end{eqnarray*}
with the choice $z=\sqrt{b\eta(r-b)}$. This yields 
\begin{eqnarray*}
B_{b,q}(N) &\prec&\frac{1}{N^{q-1}}\sum_{r=0}^{b-1}(r+1)^{q-2}+ \frac{L(b)%
\sqrt{b}}{N^{q-1}} \sum_{r=b}^{N-1} \frac{\sqrt{\eta(r-b)} }{(r+1)^{2-q}}+
r_b \\
&\prec&\left(\frac{b}{N}\right)^{q-1}\Bigl(1+\frac{L(b)}{\sqrt{b}}%
\sum_{t=0}^{N-b-1}\sqrt{\eta(t)}\Bigr) +\frac{L(b)\sqrt{b}}{N^{q-1}}%
\sum_{t=0}^{N-b-1}\frac{\sqrt{\eta(t)}}{(t+1)^{2-q}} \\
&&+\, r_b,
\end{eqnarray*}
where the second inequality follows from the change in variable $r=t+b$. We
use here $N=n$. Now if we assume that $b\prec n^{1-\delta}$, we deduce that $%
(b/n)^{q-1} \prec n^{-q\delta/2}$ and if $bL(b)\prec n^{1-\delta}$ we
analogously derive that $(b/N)^{q-1}(L(b)/\sqrt{b})\prec n^{-\frac q2
\delta} $. If $\eta(t) \prec n^{-\eta}$ and $\sigma \leq \eta/2$ we assume
that 
\begin{equation*}
\frac bn\left(1\vee \frac{L(b)}{\sqrt{b}}\right)+r_b\prec n^{-\delta},
\qquad \sum_{t=0}^{\infty}(t+1)^{q-2}\eta(t)^{1/2}<\infty.
\end{equation*}
\textbf{- }\underline{$\eta$ (\ref{nonover2},\ref{con}) case.} Note that 
\begin{eqnarray*}
C_{b,q}(r) &\prec& bL(b)\eta((r-1)b)/z+L(b)z+r_b \prec L(b) \sqrt{%
b\eta((r-1)b)}+r_b
\end{eqnarray*}
where we use $z=\sqrt{b\eta((r-1)b)}$, then 
\begin{eqnarray*}
B_{b,q}(N)&\prec& \frac{1}{N^{q-1}} +\frac{L(b)\sqrt{b}}{N^{q-1}}%
\sum_{r=1}^{N-1}\frac{\sqrt{\eta((r-1)b)}}{(r+1)^{2-q}} +r_b \\
&\prec& \frac{1}{N^{q-1}}\Bigl(1+L(b)\sqrt{b}\sum_{k=b}^{n-1}\sqrt{\eta(k)}
\Bigr)+\frac{L(b)b^{\frac 52-q}}{N^{q-1}}\sum_{k=b}^{n-1}\frac{\sqrt
{\eta(k)}}{k^{2-q}}+ r_b\!\!\! \\
\end{eqnarray*}
where the second inequality follows from replacement of $k=b(r-1)$. We use
here $N=n/b$. Let us assume $b\prec n^{1-\delta}$, then we deduce that $%
(b/n)^{q-1}\prec n^{-\frac q2\delta}$ and if $b^{3/2}L(b)\prec n^{1-\delta}$
we analogously derive that $(1/N)^{q-1}(L(b)\sqrt{b})\prec n^{-\frac
q2\delta}$. If $\eta(t) \prec n^{-\eta}$ and $\sigma \leq \eta/2$ we assume
that 
\begin{equation*}
\frac bn\left(1\vee \sqrt{b}L(b)\right)+r_b\prec n^{-\delta}, \qquad
\sum_{t=0}^{\infty}\eta(t)^{1/2}<\infty.
\end{equation*}
\textbf{- }\underline{$\lambda$ (\ref{over},\ref{con}) case.} Note that 
\begin{eqnarray*}
C_{b,q}(r)&\prec&L(b)z+\left(bL(b)/z
+\left(bL(b)/z\right)^2\right)\lambda(r-b)+r_b \\
&\prec& 2(bL(b)^2)^{\frac{2}{3}}\lambda(r-b)^{\frac{1}{3}}+(bL(b)^2)^{\frac{1%
}{3}}\lambda(r-b)^{\frac{2}{3}} + r_b,
\end{eqnarray*}
with a choice $z=(b^2L(b)\lambda(r-b))^{\frac 1{3}}$. Then 
\begin{eqnarray*}
B_{b,q}(N) &\prec&\frac{1}{N^{q-1}}\sum_{r=0}^{b-1}(r+1)^{q-2}+\frac{%
2(bL(b)^2)^{\frac{2}{3}}}{N^{q-1}}\sum_{r=b}^{N-1}\frac{\lambda(r-b)^{\frac{1%
}{3}}}{(r+1)^{2-q}} \\
&&+\,\frac{(bL(b)^2)^{\frac{1}{3}}}{N^{q-1}}\sum_{r=b}^{N-1}\frac{
\lambda(r-b)^{\frac{2}{3}}}{(r+1)^{2-q}} + r_b \\
&\prec&\left(\frac{b}{N}\right)^{q-1}\left(1+(b^{-1}L(b)^4)^{\frac{1}{3}%
}\sum_{t=0}^{N-b-1}\lambda(t)^{\frac{1}{3}} + (b^{-1}L(b))^{\frac{2}{3}%
}\sum_{t=0}^{N-b-1}\lambda(t)^{\frac{2}{3}}\right) \\
&&+\, \frac 1{N^{q-1}} \left( (bL(b)^2)^{\frac{2}{3}}\sum_{t=0}^{N-b-1}\frac{
\lambda(t)^{\frac{1}{3}}}{(t+1)^{2-q}} + (bL(b)^2)^{\frac{1}{3}%
}\sum_{t=0}^{N-b-1}\frac{ \lambda(t)^{\frac{2}{3}} }{(t+1)^{2-q}}\right) \\
&&+\,r_b
\end{eqnarray*}
where the second inequality follows from the change in variable $r=t+b$. We
use here $N=n$. Let us assume $b\prec n^{1-\delta}$, then we deduce that $%
(b/n)^{q-1} \prec n^{-\frac q2\delta}$ and if $\left( (bL(b)^{2})^{2/3} \vee
(b L(b)^2)^{1/3} \right) \prec n^{1-\delta}$ we analogously derive that $%
(b/N)^{q-1}\left((b^{-1}L(b)^4)^{\frac{1}{3}}\vee (b^{-1}L(b))^{\frac{2}{3}%
}\right) \prec n^{-\frac q2\delta}$. If $\lambda(t) \prec n^{-\lambda}$ and $%
\sigma \leq \lambda/2 $ and $\sigma \leq 2 \lambda/3$ we assume that 
\begin{equation*}
\frac bn \left(1 \vee (b^{-1}L(b)^4)^{\frac{1}{3}}\vee (b^{-1}L(b))^{\frac{2%
}{3}}\right) +r_b \prec n^{-\delta}
\end{equation*}
with $\displaystyle
\sum_{t=0}^{\infty}(t+1)^{q-2}\lambda(t)^{1/3}<\infty$ and $\displaystyle %
\sum_{t=0}^{\infty}(t+1)^{q-2}\lambda(t)^{2/3}< \infty$.\newline
\newline
\textbf{- }\underline{$\lambda$ (\ref{nonover2},\ref{con}) case.} Note that 
\begin{eqnarray*}
C_{b,q}(r)&\prec&\left(L(b)z+r_b\right)+\left(\frac{bL(b)}{z} +\left(\frac{%
bL(b)}{z}\right)^2\right)\lambda((r-1)b) \\
&\prec& 2(bL(b)^2)^{\frac{2}{3}}\lambda((r-1)b)^{\frac{1}{3}}+(bL(b)^2)^{%
\frac{1}{3}}\lambda((r-1)b)^{\frac{2}{3}}+r_b,
\end{eqnarray*}
with a choice $z=(b^2L(b)\lambda((r-1)b))^{\frac1{3}}$. Then we obtain 
\begin{eqnarray*}
B_{b,q}(N) &\prec&\frac{1}{N^{q-1}} + \frac{(bL(b)^2)^{\frac 23}}{N^{q-1}}%
\sum_{r=1}^{N-1}(r+1)^{q-2}\lambda((r-1)b)^{\frac 13} \\
&&+\, \frac{(bL(b)^2)^{\frac 13}}{N^{q-1}}\sum_{r=1}^{N-1}(r+1)^{q-2}%
\lambda((r-1)b)^{\frac 23}+ r_b \\
&\prec& \frac 1{N^{q-1}} \Bigl( 1+ (bL(b)^2)^{\frac 23 } \sum_{k=1}^{n-1}
\lambda(k)^{\frac 13}+ (bL(b)^2)^{\frac 13} \sum_{k=1}^{n-1}
\lambda(k)^{\frac 23} \Bigr) \\
&&+\, \frac 1{N^{q-1}} \Bigl( L(b)^{\frac 43}b^{\frac 83-q} \sum_{k=1}^{n-1} 
\frac{ \lambda(k)^{\frac 13}}{k^{2-q}} + L(b)^{\frac 23}b^{\frac 73-q}
\sum_{k=1}^{n-1} \frac{\lambda(k)^{\frac 23} }{k^{2-q}} \Bigr) + r_b
\end{eqnarray*}
where the second inequality follows from the change in variable $k=b(r-1)$.
We use here $N=n/b$. Let us assume $b\prec n^{1-\delta}$, then we deduce
that $(b/n)^{q-1} \prec n^{-\frac q2\delta}$ and if $\left(b^{\frac
53}L(b)^{\frac 43 } \vee b^{\frac 43}L(b)^{\frac 23}\right) \prec
n^{1-\delta}$ then we analogously derive that $(b/n)^{q-1}\left((bL(b)^2)^{%
\frac 23 } \vee (bL(b)^2)^{\frac 13}\right) \prec n^{-\frac q2\delta}$. If $%
\lambda(t) \prec n^{-\lambda}$ and $\sigma \leq \lambda/3 $ and $\sigma \leq
2 \lambda/3$ we assume that 
\begin{equation*}
\frac bn \left((bL(b)^2)^{\frac 23 } \vee (bL(b)^2)^{\frac 13}\right)
+r_b\prec n^{-\delta}
\end{equation*}
with $\displaystyle \sum_{t=0}^{n-1}(t+1)^{q-2}\lambda(t)^{1/3}<\infty$ and $%
\displaystyle \sum_{t=0}^{n-1}(t+1)^{q-2}\lambda(t)^{2/3}<\infty$ that bound
holds. \newline
Consider now the cases under assumption (\ref{conc}). \newline
\textbf{- }\underline{$\eta$ (\ref{over},\ref{conc}) case.} Note that 
\begin{eqnarray*}
C_{b,q}(r)&\prec& L(b)b\eta(r-b)/z+C(b) z^c \qquad \qquad \qquad \qquad
\qquad \qquad \qquad \qquad \\
&\prec&\Bigl(C(b)(bL(b)\eta(r-b))^c \Bigr)^{\frac
1{1+c}}+\Bigl(bL(b)\eta(r-b)\Bigr)^{\frac {2+c}{1+c}}
\end{eqnarray*}
with a choice $z=(bL(b)\eta(r-b)/C(b))^{\frac 1{c+1}}$, then 
\begin{eqnarray*}
&B_{b,q}(N)& \\
&\prec& \frac{1}{N^{q-1}}\sum_{r=0}^{b-1}(r+1)^{q-2} +\frac{%
\left(C(b)(bL(b))^c\right)^{\frac{1}{1+c}}}{N^{q-1}}%
\sum_{r=b}^{N-1}(r+1)^{q-2}(\eta(r-b))^{\frac c{1+c}} \\
&&+\frac{\left(bL(b)\right)^{\frac{2+c}{1+c}}}{N^{q-1}}%
\sum_{r=b}^{N-1}(r+1)^{ q-2} (\eta(r-b))^{\frac {2+c}{1+c}} \\
&\prec& \left(\frac{b}{N}\right)^{q-1}\Bigl(1+\left(C(b)b^{-1}L(b)^c
\right)^{\frac{1}{1+c}}\sum_{t=0}^{N-b-1}\eta(t)^{\frac{c}{1+c}} +
\left(bL(b)^{2+c}\right)^{\frac{1}{1+c}} \sum_{t=0}^{N-b-1}\eta(t)^{\frac{2+c%
}{1+c}} \Bigr) \\
&&+\, \frac{1}{N^{q-1}} \Bigl(\left(C(b)(bL(b))^c \right)^{\frac{1}{1+c}%
}\sum_{t=0}^{N-b-1}\frac{\eta(t)^{\frac{c}{1+c}}}{(t+1)^{2-q}} +
\left(bL(b)\right)^{\frac{2+c}{1+c}}\sum_{t=0}^{N-b-1}\frac{\eta(t)^{\frac{%
2+c}{1+c}} }{(t+1)^{2-q}}\Bigr) \\
\end{eqnarray*}
where the second inequality follows from the change in variable $r=t+b$. We
use here $N=n$. Now if we assume $b\prec n^{1-\delta}$ we deduce that $%
(b/n)^{q-1} \prec n^{-\frac q2 \delta}$ and if $b\left((C(b)b^{-1}L(b)^c)^{%
\frac{1}{1+c}} \vee (bL(b)^{2+c})^{\frac{2+c}{1+c}}\right) \prec
n^{1-\delta} $ we analogously derive $(b/N)^{q-1}\left((C(b)b^{-1}L(b)^c)^{%
\frac{1}{1+c}} \vee (bL(b)^{2+c})^{\frac{2+c}{1+c}}\right) \prec n^{-\frac
q2\delta}.$ \newline
If $\eta(t) \prec n^{-\eta}$ and $\sigma \leq \eta \frac c{1+c} $ and $%
\sigma \leq \eta \frac {2+c}{1+c} $ we assume that 
\begin{equation*}
\frac{b}{n}\left(1\vee\left(C(b)b^{-1}L(b)^c\right)^{\frac{1}{1+c}} \vee
\left(bL(b)^{2+c}\right)^{\frac{1}{1+c}}\right)\prec n^{-\delta}
\end{equation*}
which implies with $\displaystyle \sum_{t=0}^{\infty}\left(t+1\right)^{q-2}%
\eta(t)^{\frac{c}{1+c}}<b^{q-2}$ and $\displaystyle \sum_{t=0}^{\infty}%
\left(t+1\right)^{q-2}\eta(t)^{\frac{2+c}{1+c}}<b^{q-2}$ that $B_{b,q}(N)
\prec n^{-\frac q2\delta}.$ \newline
\textbf{- }\underline{$\eta$ (\ref{nonover2},\ref{conc}) case.} Note that 
\begin{eqnarray*}
C_{b,q}(r)&\prec& L(b)b\eta((r-1)b)/z+C(b)z^c \qquad \qquad \qquad \qquad
\qquad \qquad \qquad \qquad \qquad \qquad \qquad \qquad \\
& \prec &
(C(b)(bL(b)\eta((r-1)b))^c))^\frac1{1+c}+(bL(b)\eta((r-1)b))^{\frac
{2+c}{1+c}}
\end{eqnarray*}
with a choice $z=(bL(b)\eta((r-1)b)/C(b))^{\frac 1{c+1}}$. Then 
\begin{eqnarray*}
&B_{b,q}(N)& \\
&\prec& \frac{1}{N^{q-1}} + \frac{\left(C(b)(b L(b))^c\right)^{\frac 1{1+c}}%
}{N^{q-1}}\sum_{r=1}^{N-1}(r+1)^{q-2} \eta((r-1)b)^{\frac c{1+c}} \\
&&+\, \frac{(bL(b))^{\frac{2+c}{1+c}}}{N^{q-1}}\sum_{r=1}^{N-1}(r+1)^{q-2}
\eta((r-1)b)^{\frac {2+c}{1+c}} \\
&\prec& \frac{1}{N^{q-1}} \Bigl(1+
\left(C(b)(bL(b))^c\right)^{\frac1{1+c}}\sum_{k=b}^{n-1} \eta(k)^{\frac{c}{%
1+c}} + (bL(b))^{\frac{2+c}{1+c}}\sum_{k=b}^{n-1}\eta(k)^{\frac{2+c}{1+c}}
\Bigr) \\
&&+\, \frac{1}{N^{q-1}} \Bigl(\left(C(b) L(b)^c \right)^{\frac{1}{1+c}} b^{%
\frac{2+3c}{1+c}-q} \sum_{k=b}^{n-1}\frac{\eta(k)^{\frac{c}{1+c}}}{ k^{2-q}}
\\
&&+\, L(b)^\frac{2+c}{1+c} b^{\frac{4+3c}{1+c}-q} \sum_{k=b}^{n-1} \frac{%
\eta(k)^{\frac{2+c}{1+c}}}{k^{2-q}} \Bigr)
\end{eqnarray*}
where the second inequality folllows from the change in variables $k=b(r-1)$%
. We use here $N=n/b$. Now if we assume $b\prec n^{1-\delta}$, we deduce
that $(b/n)^{q-1} \prec n^{-\frac q2 \delta}$. Now if $\left(1\vee
(C(b)(bL(b))^c)^{\frac1{1+c}} \vee (bL(b))^{\frac{2+c}{1+c}}\right)\prec
n^{-\delta} $ we analogously derive $(1/N)^{q-1} \left(1\vee
(C(b)(bL(b))^c)^{\frac1{1+c}} \vee (bL(b))^{\frac{2+c}{1+c}}\right)\prec
n^{-\frac q2 \delta}$. If $\eta(t) \prec n^{-\eta}$ and $\sigma \leq \eta
\frac c{1+c} $ and $\sigma \leq \eta \frac {2+c}{1+c} $ we assume that 
\begin{equation*}
\frac{b}{n}\left(1\vee (C(b)(bL(b))^c)^{\frac1{1+c}} \vee (bL(b))^{\frac{2+c%
}{1+c}}\right)\prec n^{-\delta}
\end{equation*}
which implies with $\displaystyle
\sum_{t=0}^{n-1}\left(t+1\right)^{q-2}\eta(t)^{\frac{c}{1+c}}\prec b^{q-2}$
and $\displaystyle \sum_{t=0}^{n-1}\left(t+1\right)^{q-2}\eta(t)^{\frac{2+c}{%
1+c}}\prec b^{q-2}$ that $B_{b,q}(N) \prec n^{-\frac q2\delta}.$ \newline
\textbf{- }\underline{$\lambda$ (\ref{over},\ref{conc}) case.} Note that 
\begin{eqnarray*}
C_{b,q}(r)&\prec& C(b)z^c+\left(bL(b)/z+\left(bL(b)/z
\right)^2\right)\lambda(r-b), \\
&\prec& (C(b)bL(b)^c )^{\frac{2}{2+c}}(\lambda(r-b))^{\frac{c}{2+c}} +
\left(C(b)(bL(b)^c)^{\frac{1}{2+c}}\right)(\lambda(r-b))^{\frac{1+c}{2+c}} \\
&&+\, \left(C(b)(bL(b)^{2c})^{\frac{1}{2+c}}\right)(\lambda(r-b))^{\frac{1+c%
}{2+c}},
\end{eqnarray*}
with a choice $z=((bL(b))^2 C(b)^{-1}\lambda(r-b))^{\frac1{2+c}}$. 
\begin{eqnarray*}
&B_{b,q}(N)& \\
&\prec& \frac {1}{N^{q-1}}\sum_{r=0}^{b-1}(r+1)^{q-2}+\frac{(C(b)(bL(b))^c)^{%
\frac{2}{2+c}}}{N^{q-1}}\sum_{r=b}^{N-1}\frac{ \lambda(r-b)^{\frac{c}{2+c}}}{%
(r+1)^{2-q}} \\
&&+\, \frac{(C(b)(bL(b))^c)^{\frac{1}{2+c}}}{N^{q-1}}\sum_{r=b}^{N-1}\frac{
\lambda(r-b)^{\frac{1+c}{2+c}}}{(r+1)^{2-q}} + \frac{(C(b)(bL(b))^{2c})^{%
\frac{1}{2+c}}}{N^{q-1}}\sum_{r=b}^{N-1} \frac{ \lambda(r-b)^{\frac{1+c}{2+c}%
} }{(r+1)^{2-q}} \\
&\prec& \left(\frac{b}{N}\right)^{q-1} \Bigr(1+(C(b)^2b^{c-2}L(b)^{2c})^{%
\frac{1}{2+c}} \sum_{t=0}^{N-b-1}\lambda(t)^{\frac{c}{2+c}} \\
&&+\,(C(b)b^{-2}L(b)^c)^{\frac{1}{2+c}}\sum_{t=0}^{N-b-1}\lambda(t)^{\frac{%
1+c}{2+c}} + (C(b)b^{c-2}L(b)^{2c})^{\frac{1}{2+c}}\sum_{t=0}^{N-b-1}%
\lambda(t)^{\frac{1+c}{2+c}} \Bigl) \\
&&+\, \frac{1}{N^{q-1}}\Bigl( (C(b)(bL(b))^c)^{\frac{2}{2+c}%
}\sum_{t=0}^{N-b-1} \frac{ \lambda(t)^{\frac{c}{2+c}} }{(t+1)^{2-q}}
+(C(b)(bL(b))^c)^{\frac{1}{2+c}}\sum_{t=0}^{N-b-1} \frac{ \lambda(t)^{\frac{%
1+c}{2+c}} }{(t+1)^{2-q}} \\
&&+\, (C(b)(bL(b))^{2c})^{\frac{1}{2+c}}\sum_{t=0}^{N-b-1} \frac{
\lambda(t)^{\frac{1+c}{2+c}} }{(t+1)^{2-q}} \Bigr)
\end{eqnarray*}
where the second inequality follows from the change in variable $r=t+b$. We
use here $N=n$. Now if we assume that $b\prec n^{1-\delta}$ we deduce that $%
(b/n)^{q-1} \prec n^{-\frac q2 \delta}$.\newline
If $b\left((C(b)^2b^{c-2}L(b)^{2c})^{\frac{1}{2+c}}\vee(C(b)b^{-2}L(b)^c)^{%
\frac{1}{2+c}}\vee(C(b)b^{c-2}L(b)^{2c})^{\frac{1}{2+c}}\right) \prec
n^{1-\delta}$ we analogously derive 
\begin{equation*}
(b/N)^{q-1}\left((C(b)^2b^{c-2}L(b)^{2c})^{\frac{1}{2+c}%
}\vee(C(b)b^{-2}L(b)^c)^{\frac{1}{2+c}}\vee(C(b)b^{c-2}L(b)^{2c})^{\frac{1}{%
2+c}}\right) \prec n^{-\frac q2\delta}.
\end{equation*}
If $\lambda(t) \prec n^{-\lambda}$ and $\sigma \leq \lambda \frac c{1+c} $
and $\sigma \leq \lambda \frac {2+c}{1+c} $ we assume that 
\begin{equation*}
\frac bn\left(1 \vee (C(b)^2b^{c-2}L(b)^{2c})^{\frac{1}{2+c}%
}\vee(C(b)b^{-2}L(b)^c)^{\frac{1}{2+c}}\vee(C(b)b^{c-2}L(b)^{2c})^{\frac{1}{%
2+c}} \right)\prec n^{-\delta}
\end{equation*}
which implies with $\displaystyle
\sum_{t=0}^{\infty}\left(t+1\right)^{q-2}\lambda(t)^{\frac{c}{1+c}}<b^{q-2}$
and $\displaystyle
\sum_{t=0}^{\infty}\left(t+1\right)^{q-2}\lambda(t)^{\frac{2+c}{1+c}%
}<b^{q-2} $ that $B_{b,q}(N) \prec n^{-\frac q2\delta}.$\newline
\textbf{- }\underline{$\lambda$ (\ref{nonover2},\ref{conc}) case.} Note that 
\begin{eqnarray*}
C_{b,q}(r)&\prec& C(b)z^c+\left(\frac{bL(b)}{z}+\left(\frac{bL(b)}{z}%
\right)^2\right)\lambda((r-1)b) \\
&\prec& (C(b)(bL(b))^c)^{\frac{2}{2+c}}\lambda((r-1)b)^{\frac{c}{2+c}} +
\left((C(b)(bL(b))^c)^{\frac{1}{2+c}}\right)\lambda((r-1)b)^{\frac{1+c}{2+c}}
\\
&&+\, \left((C(b)(bL(b))^{2c})^{\frac{1}{2+c}}\right)\lambda((r-1)b)^{\frac{%
1+c}{2+c}},
\end{eqnarray*}
with a choice $z=((bL(b))^2C(b)^{-1}\lambda((r-1)b))^{\frac 1{2+c}}$. We
obtain 
\begin{eqnarray*}
B_{b,q}(N)&\prec& \frac {1}{N^{q-1}}+ \frac{(C(b)(bL(b))^c)^{\frac{2}{2+c}}}{%
N^{q-1}}\sum_{r=1}^{N-1}(r+1)^{q-2}\lambda((r-1)b)^{\frac{c}{2+c}} \\
&&+\, \frac{(C(b)(bL(b))^c)^{\frac{1}{2+c}}}{N^{q-1}}%
\sum_{r=1}^{N-1}(r+1)^{q-2}\lambda((r-1)b)^{\frac{1+c}{2+c}} \\
&&+\,\frac{(C(b)(bL(b))^{2c})^{\frac{1}{2+c}}}{N^{q-1}}
\sum_{r=1}^{N-1}(r+1)^{q-2}\lambda((r-1)b)^{\frac{1+c}{2+c}} \\
\\
&\prec& \frac{1}{N^{q-1}}\Bigr(1+(C(b)(bL(b))^c)^{\frac{2}{2+c}%
}\sum_{k=1}^{n-1}\lambda(k)^{\frac{c}{2+c}} + (C(b)(bL(b))^c)^{\frac{1}{2+c}%
}\sum_{k=1}^{n-1}\lambda(k)^{\frac{1+c}{2+c}} \\
&&+\, (C(b)(bL(b))^{2c})^{\frac{1}{2+c}}\sum_{k=1}^{n-1}\lambda(k)^{\frac{1+c%
}{2+c}} \Bigl) + \frac{1}{N^{q-1}}\Bigr((C(b)L(b)^c)^{\frac{2}{2+c}}b^{\frac{%
4(1+c)}{2+c}-q} \sum_{k=1}^{n-1}\frac{\lambda(k)^{\frac{c}{2+c}}}{k^{2-q} }
\\
&&+\, (C(b)L(b)^c)^{\frac{1}{2+c}}b^{\frac{4+3c}{2+c}-q} \sum_{k=1}^{n-1}%
\frac{\lambda(k)^{\frac{1+c}{2+c}}}{k^{2-q} } + (C(b)L(b)^{2c})^{\frac{1}{2+c%
}}b^{\frac{4(1+c)}{2+c}-q}\sum_{k=1}^{n-1}\frac{ \lambda(k)^{\frac{1+c}{2+c}}%
}{k^{2-q} } \Bigr)
\end{eqnarray*}
where the second inequality follows from the change in variable $k=b(r-1)$.
We use here $N=n/b$. Now if we assume that $b\prec n^{1-\delta}$ we deduce
that $(b/n)^{q-1} \prec n^{-\frac q2 \delta}$.\newline
If $b\left((C(b)(bL(b))^c)^{\frac{2}{2+c}}\vee (C(b)(bL(b))^c)^{\frac{1}{2+c}%
}\vee \vee (C(b)(bL(b))^{2c})^{\frac{1}{2+c}} \right) \prec n^{1-\delta}$ we
analogously derive that \newline
$(b/N)^{q-1}\left((C(b)(bL(b))^c)^{\frac{2}{2+c}}\vee (C(b)(bL(b))^c)^{\frac{%
1}{2+c}}\vee (C(b)(bL(b))^{2c})^{\frac{1}{2+c}} \right) \prec n^{-\frac q2
\delta}$. If $\lambda(t) \prec n^{-\lambda}$ and $\sigma \leq \lambda \frac
c{1+c} $ and $\sigma \leq \lambda \frac {2+c}{1+c} $ we assume that 
\begin{equation*}
\frac bn \left(1 \vee (C(b)(bL(b))^c)^{\frac{2}{2+c}}\vee (C(b)(bL(b))^c)^{%
\frac{1}{2+c}}\vee (C(b)(bL(b))^{2c})^{\frac{1}{2+c}} \right) \prec n^{-
\delta}
\end{equation*}
which implies with $\displaystyle
\sum_{t=0}^{n-1}(t+1)^{q-2}\lambda(t)^{\frac{c}{1+c}}\prec b^{q-2}$ and $%
\displaystyle
\sum_{t=0}^{n-1}(t+1)^{q-2}\lambda(t)^{\frac{2+c}{1+c}}\prec b^{q-2}$ that $%
B_{b,q}(N) \prec n^{-\frac q2\delta}.$

\begin{Lemma}
The relation $\widehat{\Delta }_{b,n}^{(2)}(x)\rightarrow _{n\rightarrow
\infty }0$ holds in the following cases under the convergence assumption (%
\ref{con})

\begin{itemize}
\item  In the overlapping case, if we have respectively 
\begin{eqnarray*}
&\sum_{t=0}^{\infty }\eta (t)^{1/2}<\infty ,\mbox{and}&r_{b}+\frac{b}{n}%
\left( 1\vee \frac{L(b)}{\sqrt{b}}\right)\to0 ,\qquad \qquad \quad , \\
&\sum_{t=0}^{\infty }\lambda (t)^{2/3}<\infty ,\mbox{and}&r_{b}+\frac{b}{n}%
\left( 1\vee \left( \frac{L(b)^{4}}{b}\right) ^{1/3}\vee \left( \frac{L(b)}{b%
}\right) ^{2/3}\right) \to0 .
\end{eqnarray*}

\item  In the non-overlapping case, if we have respectively 
\begin{eqnarray*}
\sum_{t=0}^{\infty }\eta (t)^{1/2}< \infty ,&\mbox{ and }&r_{b}+\frac{b}{n}%
\left( 1\vee \sqrt{b}L(b)\right)\to0 , \\
\sum_{t=0}^{\infty }\lambda (t)^{2/3}< \infty ,&\mbox{and }&r_{b}+\frac{b}{n}%
\left( 1\vee (bL(b)^{2})^{2/3}\vee (bL(b)^{2})^{1/3}\right)\to0 .
\end{eqnarray*}
\end{itemize}
\end{Lemma}

This lemma together with lemma \ref{gen_b} yields theorem \ref{rsubsamp2}.

\begin{Lemma}
The relation $B_{b,q}(N)\prec n^{-q\delta /2}$ holds under concentration
assumption (\ref{conc}) if respectively \emph{the overlapping setting is used%
} and one among the following relations hold as $n\rightarrow \infty $ 
\begin{equation*}
\eta \mbox{-dependence:}\sum_{t=0}^{\infty }(t+1)^{q-2}\eta (t)^{\frac{2+c}{%
1+c}}<\infty ,\qquad \qquad \qquad \qquad \qquad \qquad \qquad \qquad \qquad
\end{equation*}
\begin{equation*}
\frac{b}{n}\left( 1\vee \left( C(b)b^{-1}L(b)^{c}\right) ^{\frac{1}{1+c}%
}\vee \left( bL(b)^{2+c}\right) ^{\frac{1}{1+c}}\right) \prec n^{-\delta },
\end{equation*}
\begin{equation*}
\lambda \mbox{-dependence:}\sum_{t=0}^{\infty }(t+1)^{q-2}\lambda (t)^{\frac{%
1+c}{2+c}}<\infty ,\mbox{ }\qquad \qquad \qquad \qquad \qquad \qquad \qquad
\qquad \qquad
\end{equation*}
\begin{equation*}
\frac bn\left(1 \vee (C(b)^2b^{c-2}L(b)^{2c})^{\frac{1}{2+c}%
}\vee(C(b)b^{-2}L(b)^c)^{\frac{1}{2+c}}\vee(C(b)b^{c-2}L(b)^{2c})^{\frac{1}{%
2+c}} \right)\prec n^{-\delta}
\end{equation*}
or \emph{the non-overlapping setting is used} and 
\begin{equation*}
\eta \mbox{-dependence:}\sum_{t=0}^{n-1}(t+1)^{q-2}\eta (t)^{\frac{2+c}{1+c}%
}\prec b^{q-2} ,\qquad \qquad \qquad \qquad \qquad \qquad \qquad \qquad
\qquad
\end{equation*}
\begin{equation*}
\frac{b}{n}\left( 1\vee \left( C(b)(bL(b))^{c}\right) ^{\frac{1}{1+c}}\vee
(bL(b))^{\frac{2+c}{1+c}}\right) \prec n^{-\delta },
\end{equation*}
\begin{equation*}
\lambda \mbox{-dependence:}\sum_{t=0}^{n-1}(t+1)^{q-2}\lambda (t)^{\frac{1+c%
}{2+c}}\prec b^{q-2} ,\qquad \qquad \qquad \qquad \qquad \qquad \qquad
\qquad \qquad
\end{equation*}
\begin{equation*}
\frac{b}{n}\left( 1\vee (C(b)(bL(b))^{c})^{\frac{2}{2+c}}\vee
(C(b)(bL(b))^{c})^{\frac{1}{2+c}}\vee (C(b)(bL(b))^{2c})^{\frac{1}{2+c}%
}\right) \prec n^{-\delta }.
\end{equation*}
\end{Lemma}

This lemma together with lemma \ref{gen_b} yields theorem \ref{rsubsamp}.


\subsection{Proof of Theorem \ref{deftheta}}

Put $k_{n}=\left\lfloor n/a_{n}\right\rfloor $. Partition $\{1,\ldots ,n\}$
into $k_{n}$ blocks of size $a_{n}$%
\begin{equation*}
J_{j}=J_{j,n}=\left\{ \left( j-1\right) a_{n}+1,\ldots,ja_{n}\right\}
,\qquad j=1,\ldots,k_{n},
\end{equation*}
and, in case $k_{n}a_{n}<n$, a remainder block, $J_{k_{n}+1}=\left\{
k_{n}a_{n}+1,\ldots,n\right\} $. Observe that 
\begin{equation*}
\mathbb{P}(M_{n}\leq w_{n}(x))=\mathbb{P}\left(
\bigcap_{j=1}^{k_{n}+1}\left\{ M\left( J_{j}\right) \leq w_{n}(x)\right\}
\right)
\end{equation*}
where $M\left( J_{j}\right) =\max_{i\in J_{j}}(X_{i})$. Since $\mathbb{P}%
(M\left( J_{j}\right) >w_{n}(x))\leq a_{n}\bar{F}\left( w_{n}(x)\right)
\rightarrow 0$ as $n\rightarrow \infty $, the remainder block can be omitted
and 
\begin{equation*}
\mathbb{P}(M_{n}\leq w_{n}(x))=\mathbb{P}\left( \bigcap_{j=1}^{k_{n}}\left\{
M\left( J_{j}\right) \leq w_{n}(x)\right\} \right) +o\left( 1\right) .
\end{equation*}
Let 
\begin{eqnarray*}
J_{j}^{\ast } &=&J_{j,n}^{\ast }=\left\{ \left( j-1\right)
a_{n}+1,\ldots,ja_{n}-l_{n}\right\} ,\qquad j=1,\ldots,k_{n}, \\
J_{j}^{\prime } &=&J_{j,n}^{\prime }=\left\{
ja_{n}-l_{n},\ldots,ja_{n}\right\} ,\qquad j=1,\ldots,k_{n}.
\end{eqnarray*}
Since $P\left( \bigcup_{j=1}^{k_{n}}M\left( J_{j}^{\prime }\right)
>w_{n}(x)\right) \leq k_{n}l_{n}\bar{F}\left( w_{n}(x)\right) \rightarrow 0$
as $n\rightarrow \infty $, we deduce that 
\begin{equation*}
\mathbb{P}(M_{n}\leq w_{n}(x))=\mathbb{P}\left( \bigcap_{j=1}^{k_{n}}\left\{
M\left( J_{j}^{\ast }\right) \leq w_{n}(x)\right\} \right) +o\left( 1\right)
.
\end{equation*}
Let $B_{j}=B_{j,n}=\left\{ M\left( J_{j}^{\ast }\right) \leq
w_{n}(x)\right\} $. We write 
\begin{eqnarray*}
&&\mathbb{P}\left( \bigcap_{j=1}^{k_{n}}B_{j}\right) -\prod_{j=1}^{k_{n}}%
\mathbb{P}(B_{j}) \\
&=&\sum_{i=1}^{k_{n}}\left( \mathbb{P}\left(
\bigcap_{j=1}^{k_{n}-i+1}B_{j}\right) \prod_{j=k_{n}-i+2}^{k_{n}}\mathbb{P}%
(B_{j})-\mathbb{P}\left( \bigcap_{j=1}^{k_{n}-i}B_{j}\right)
\prod_{j=k_{n}-i+1}^{k_{n}}\mathbb{P}(B_{j})\right) \\
&=&\sum_{i=1}^{k_{n}}\left( \mathbb{P}\left(
\bigcap_{j=1}^{k_{n}-i+1}B_{j}\right) -\mathbb{P}\left(
\bigcap_{j=1}^{k_{n}-i}B_{j}\right) \mathbb{P}(B_{k_{n}-i+1})\right)
\prod_{j=k_{n}-i+2}^{k_{n}}\mathbb{P}(B_{j}).
\end{eqnarray*}
We want to bound the following quantity 
\begin{equation*}
\left| \mathbb{P}\left( \bigcap_{j=1}^{k_{n}-i+1}B_{j}\right) -\mathbb{P}%
\left( \bigcap_{j=1}^{k_{n}-i}B_{j}\right) \mathbb{P}(B_{k_{n}-i+1})\right| .
\end{equation*}
Let us define $f_{n}^{\left( x\right) }(y)=\mathbb{I}_{\{y\leq w_{n}(x)\}}$.
Let $(\alpha _{n})$ be a sequence such that $\alpha _{n}\rightarrow 0$ as $%
n\rightarrow \infty $ and put $x_{n}^{-}=x-\alpha _{n}$ and $%
x_{n}^{+}=x+\alpha _{n}$. We simply approximate the function $f_{n}^{\left(
x\right) }$ by Lipschitz and bounded functions $g_{n},h_{n}\in \mathcal{F}%
_{1}$ with 
\begin{equation*}
f_{n}^{(x_{n}^{-})}\leq g_{n}\leq f_{n}^{\left( x\right) }\leq h_{n}\leq
f_{n}^{(x_{n}^{+})}
\end{equation*}
and we quote that it is easy to choose functions $g_{n}$ and $h_{n}$ with
Lipschitz coefficient $u_{n}\alpha _{n}^{-1}$. For $I\subset \{1,\ldots ,n\}$%
, let $H_{I}(f_{n}^{\left( x\right) })=\mathbb{E}\left[ \prod_{i\in
I}f_{n}^{\left( x\right) }(X_{i})\right] $. Note that 
\begin{equation*}
H_{I}(f_{n}^{\left( x_{n}^{-}\right) })\leq H_{I}(g_{n})\leq
H_{I}(f_{n}^{\left( x\right) })\leq H_{I}(h_{n})\leq H_{I}(f_{n}^{\left(
x_{n}^{+}\right) }).
\end{equation*}
Let $C_{I,J}(f_{n}^{\left( x\right) })=H_{I\cup J}(f_{n}^{\left( x\right)
})-H_{I}(f_{n}^{\left( x\right) })H_{J}(f_{n}^{\left( x\right) })$, we have 
\begin{equation*}
C_{I,J}(g_{n})-\delta _{I,J}\left( g_{n},h_{n}\right) \leq
C_{I,J}(f_{n}^{\left( x\right) })\leq C_{I,J}(h_{n})+\delta _{I,J}\left(
g_{n},h_{n}\right)
\end{equation*}
with 
\begin{equation*}
\delta _{I,J}\left( g_{n},h_{n}\right)
=H_{I}(h_{n})H_{J}(h_{n})-H_{I}(g_{n})H_{J}(g_{n}).
\end{equation*}

Let $I_{i}=\left\{ l:\{X_{l}\leq w_{n}(x)\}\in \bigcap
_{j=1}^{k_{n}-i}B_{j}\right\} $ and $J_{i}=\left\{ l:\{X_{l}\leq
w_{n}(x)\}\in B_{k_{n}-i+1}\right\} $. We have 
\begin{eqnarray*}
\left| H_{I_{i}}(h_{n})-H_{I_{i}}(g_{n})\right| &\leq &\left(
k_{n}-i+1\right) a_{n}\left( \bar{F}(w_{n}(x_{n}^{-}))-\bar{F}%
(w_{n}(x_{n}^{+}))\right) \\
\left| H_{J_{i}}(h_{n})-H_{J_{i}}(g_{n})\right| &\leq &a_{n}\left( \bar{F}%
(w_{n}(x_{n}^{-}))-\bar{F}(w_{n}(x_{n}^{+}))\right)
\end{eqnarray*}
Then we have 
\begin{equation*}
\left| C_{I_{i},J_{i}}(f_{n}^{\left( x\right) })\right| \leq \left|
C_{I_{i},J_{i}}(h_{n})\right| \vee \left| C_{I_{i},J_{i}}(g_{n})\right|
+|\delta _{I_{i},J_{i}}\left( g_{n},h_{n}\right) |
\end{equation*}
and 
\begin{equation*}
|\delta _{I_{i},J_{i}}\left( g_{n},h_{n}\right) |\leq \left|
H_{I_{i}}(h_{n})-H_{I_{i}}(g_{n})\right| +\left|
H_{J_{i}}(h_{n})-H_{J_{i}}(g_{n})\right| .
\end{equation*}
Note that as $n\rightarrow \infty $%
\begin{equation*}
n\left( \bar{F}(w_{n}(x_{n}^{-}))-\bar{F}(w_{n}(x_{n}^{+}))\right) \sim
2\alpha _{n}\gamma (1+\gamma x)_{+}^{-1/\gamma -1}.
\end{equation*}
If $X$ is $\eta $-weakly dependent, it follows that 
\begin{equation*}
\left| C_{I_{i},J_{i}}(f_{n}^{\left( x\right) })\right| \leq \left(
k_{n}-i+2\right) a_{n}u_{n}\alpha _{n}^{-1}\eta (l_{n})+2\alpha _{n}\gamma
(1+\gamma x)_{+}^{-1/\gamma -1}\frac{\left( k_{n}-i+2\right) a_{n}}{n}.
\end{equation*}
An optimal choice of $\alpha _{n}$ is then given by 
\begin{equation*}
\alpha _{n}\sim \left[ n\eta (l_{n})u_{n}\right] ^{1/2}
\end{equation*}
and then 
\begin{equation*}
\left| C_{I_{i},J_{i}}(f_{n}^{\left( x\right) })\right| \prec \left( n\eta
(l_{n})u_{n}\right) ^{1/2}.
\end{equation*}
It follows that 
\begin{equation*}
\left| \mathbb{P}\left( \bigcap_{j=1}^{k_{n}}B_{j}\right)
-\prod_{j=1}^{k_{n}}\mathbb{P}(B_{j})\right| \prec k_{n}\left( n\eta
(l_{n})u_{n}\right) ^{1/2}.
\end{equation*}
If $X$ is $\lambda $-weakly dependent, it follows that 
\begin{eqnarray*}
\left| C_{I_{i},J_{i}}(f_{n}^{\left( x\right) })\right| &\leq & \left[
\left( k_{n}-i+2\right) a_{n}u_{n}\alpha _{n}^{-1}+\left( k_{n}-i+1\right)
a_{n}u_{n}^{2}\alpha _{n}^{-2}\right] \lambda (l_{n}) \\
&&+2\alpha _{n}\gamma (1+\gamma x)_{+}^{-1/\gamma -1}\frac{\left(
k_{n}-i+2\right) a_{n}}{n}.
\end{eqnarray*}
An optimal choice of $\alpha _{n}$ is then given by 
\begin{equation*}
\alpha _{n}\sim \left[ k_{n}\lambda (l_{n})u_{n}\right] ^{1/2}\vee \left[
k_{n}\lambda (l_{n})u_{n}^{2}\right] ^{1/3}
\end{equation*}
and then 
\begin{equation*}
\left| C_{I_{i},J_{i}}(f_{n}^{\left( x\right) })\right| \prec (\left[
n\lambda (l_{n})u_{n}\right] ^{1/2}\vee \left[ na_{n}\lambda (l_{n})u_{n}^{2}%
\right] ^{1/3}).
\end{equation*}
It follows that 
\begin{equation*}
\left| \mathbb{P}\left( \bigcap_{j=1}^{k_{n}}B_{j}\right)
-\prod_{j=1}^{k_{n}}\mathbb{P}(B_{j})\right| \prec k_{n}(\left[ n\lambda
(l_{n})u_{n}\right] ^{1/2}\vee \left[ na_{n}\lambda (l_{n})u_{n}^{2}\right]
^{1/3}).
\end{equation*}

Finally we deduce that 
\begin{equation*}
\mathbb{P}(M_{n}\leq w_{n}\left( x\right) )=\left[ \mathbb{P}(M_{n}\leq
w_{a_{n}}\left( x\right) )\right] ^{k_{n}}+o(1)
\end{equation*}
and the result follows.

\paragraph{Aknowledgements.}

We want to thank Patrice Bertail for essential discussions during many years
at ENSAE. The first author also wants to thank Paul Embrechts at ETHZ and
the Swiss Banking Institute at Zurich University for their strong support.

{\small \strut }

{\small Paul DOUKHAN (doukhan@u-cergy.fr) }

{\small University Cergy Pontoise, UFR Sciences-Techniques, Saint-Martin }

{\small 2, avenue Adolphe-Chauvin, B.P. 222 Pontoise }

{\small 95302 Cergy-Pontoise cedex }

{\small France }

{\small \strut }

{\small Silika PROHL (prohl@isb.uzh.ch) }

{\small Swiss Institute for International Economics and Applied Economic
Research }

{\small Zurich University }

{\small Switzerland }

{\small \strut }

{\small Christian Y. ROBERT (chrobert@ensae.fr) }

{\small Ecole Nationale de la Statistique et de l'Administration Economique }

{\small Timbre J120, 3 Avenue Pierre Larousse }

{\small 92245 Malakoff Cedex }

{\small France }

\end{document}